\documentclass{amsart}

\usepackage{amsfonts}
\usepackage{graphicx}
\usepackage{epstopdf}
\usepackage{algorithmic}
\usepackage{longtable}

\usepackage{amsopn}
\usepackage{amsmath}
\usepackage{amsfonts}
\usepackage{amssymb}
\usepackage{amscd}
\usepackage{mathtools}

\usepackage{booktabs}
\usepackage{color}
\usepackage{float}

\ifpdf
\DeclareGraphicsExtensions{.eps,.pdf,.png,.jpg}
\else
\DeclareGraphicsExtensions{.eps}
\fi

\newcommand{\reff}[1]{(\ref{#1})}
\usepackage{cleveref}
\newtheorem{theorem}{Theorem}[section]
\theoremstyle{plain}
\newtheorem{thm}[theorem]{Theorem}

\newtheorem{alg}[theorem]{Algorithm}

\newtheorem{defi}[theorem]{Definition}

\newtheorem{eg}[theorem]{Example}

\newcommand{\re}{\mathbb{R}}

\newcommand{\N}{\mathbb{N}}

\newcommand{\lmd}{\lambda}

\newcommand{\eps}{\epsilon}

\newcommand{\dt}{\delta}

\newcommand{\st}{\mathit{s.t.}}
\def\af{\alpha}

\newcommand{\mc}[1]{\mathcal{#1}}

\newcommand{\be}{\begin{equation}}
\newcommand{\ee}{\end{equation}}
\newcommand{\beq}{\begin{equation}}
\newcommand{\eeq}{\end{equation}}
\newcommand{\baray}{\begin{array}}
\newcommand{\earay}{\end{array}}

\newcommand{\bbm}{\begin{bmatrix}}
\newcommand{\ebm}{\end{bmatrix}}

\newcommand{\bit}{\begin{itemize}}
\newcommand{\eit}{\end{itemize}}

\newcommand{\bdes}{\begin{description}}
\newcommand{\edes}{\end{description}}

\numberwithin{equation}{section}

\begin{document}

\title[Polynomial Optimization Relaxations for GSIPs]
{Polynomial Optimization Relaxations for Generalized Semi-Infinite Programs}

\author{Xiaomeng Hu and Jiawang Nie}

\address{Xiaomeng Hu and Jiawang Nie,
Department of Mathematical Sciences,
University of California San Diego,
9500 Gilman Drive, La Jolla, CA 92093, USA.
}

\email{x8hu@ucsd.edu, njw@math.ucsd.edu}

	

\begin{abstract}
This paper studies generalized semi-infinite programs (GSIPs)
given by polynomials. We propose a hierarchy of
polynomial optimization relaxations to solve them.
They are based on Lagrange multiplier expressions and polynomial extensions.
Moment-SOS relaxations are applied to solve the polynomial optimization.
The convergence of this hierarchy is shown under certain conditions.
In particular, the classical semi-infinite programs (SIPs)
can be solved as a special case of GSIPs.
We also study GSIPs that have convex infinity constraints
and show that they can be solved exactly by
a single polynomial optimization relaxation.
The computational efficiency is demonstrated
by extensive numerical results.
\end{abstract}

\maketitle

\section{Introduction}
\label{sec:Intro}

We consider the generalized semi-infinite program (GSIP):
\be \label{primal-GSIP}
\left\{
\begin{array}{cl}
	\min \, & f(x) \\
	\st  \, &  g(x,u)\geq 0\quad\forall u\in U(x),\\
	&  x\in X,\\
\end{array}
\right.
\ee
where the constraining sets $X$ and $U(x)$ are given as
\be  \label{set:X}
X=\left\{ x\in\mathbb{R}^n \;\middle\vert\;
\begin{array}{lc}
	c_i(x)= 0 \,(i\in \mathcal{I}_1),\\
	c_j(x)\geq 0 \,(j\in\mathcal{J}_1)
\end{array}\right\} ,
\ee
\be \label{U(x)set}
U(x)=\left\{ u\in\mathbb{R}^p \;\middle\vert\;
\begin{array}{l}
	h_i(x,u)= 0\, (i\in \mathcal{I}_2),\\
	h_j(x,u)\geq 0\, (j\in\mathcal{J}_2)
\end{array}
\right\} .
\ee
The constraining function $g(x,u)$ in \reff{primal-GSIP}
is the $s$-dimensional vector
\[
g(x,u) \, \coloneqq \, \big( g_1(x,u), \ldots,  g_s(x,u) \big) .
\]
In this paper, we assume that all functions are given by polynomials, and
$\mathcal{I}_1,\,\mathcal{I}_2,\,\mathcal{J}_1$ and $\mathcal{J}_2$
are finite label sets.
For each feasible $x$, the constraint $g(x,u) \ge 0$
is required to hold for {\it every} point $u \in U(x)$.
In particular, for the special case $U(x) = \emptyset$ (the empty set),
this constraint automatically holds for all $x \in X$.
Throughout the paper, $a \ge b$ means that
the entries of $a-b$ are all nonnegative.
Generally, $U(x)$ is an infinite set,
which makes it challenging to solve the GSIP.
The feasible set of \eqref{primal-GSIP} is denoted as
\be \label{gsipset:F}
\mathcal{F} \, \coloneqq \,
\big \{  x \in X :  g(x,u)\geq 0\,\forall u\in U(x) \big \} .
\ee
The problem~\eqref{primal-GSIP} is called a
semi-infinite program (SIP) if the set $U(x)\equiv U$
is independent of $x$.

GSIPs have broad applications.
One is the centering design problem,
which aims to maximize the volume of a parameterized set $U(x)$,
such as  balls/ellipsoids \cite{boyd2004convex},
boxes \cite{bemporad2004inner} and gemstones \cite{winterfeld2008cutting},
inscribed in a given region \cite{stein2002GSIPBilevel}.
For instance, an ellipsoid can be expressed as
\[
U(x) \, = \, \big \{ u \in \re^p : \,
( u - a(x) )^T P(x)^{-1}( u - a(x) )\leq 1
\big \},
\]
where $a(x)$ is a vector parameterizing the center and
$P(x)$ is a symmetric positive definite matrix
(i.e., $P(x) \succ 0$ for $x \in X$)
parameterizing the ellipsoid shape.
The maximum ellipsoid design problem is to
look for $x \in X$ such that $U(x)$ has maximum volume and
$U(x) \subseteq Y$, for a given set $Y \subseteq \re^p$.
If $Y = \{u: g(u) \ge 0 \}$ for a vector function $g$,
this design problem can be formulated as the GSIP
(see \cite{stein2002GSIPBilevel,stein2003solvingpaper}):
\be  \label{design centering problem}
\left\{
\begin{array}{cl}
	\max \, & \det \, P(x) \\
	\st  \, & g(u) \ge 0 \,\, \forall \, u \in U(x), \\
	& x \in X.
\end{array}
\right.
\ee
Other applications include
min-max optimization \cite{pang2020discretization,zakovic2003semi},
portfolio optimization \cite{lemonidis2008global},
disjunctive programming \cite{kirst2016disjunctive},
the Stackelberg game \cite{stein2002GSIPBilevel},
and dynamical network classification\cite{weber2007genegeneralized}.
More applications can be found in
\cite{vazquez2008generalized,spissueLQW08,still1999generalized}.

There are various methods for solving SIPs
(i.e., $U(x) \equiv U$ is constant).
Some methods first approximate $U$ by a finite subset and
then solve the resulting approximation optimization
problems with finitely many constraints.
They subsequently update the approximation
until a certain termination criterion is met \cite{kirst2016disjunctive}.
When SIPs are given by polynomials, a semidefinite relaxation
method is given in \cite{wangliguofeng2014semidefinite}.
Different methods can be obtained
when $U$ is approximated by different schemes.
Typical tricks for the approximation are the exchange,
discretization and local reduction. We refer to \cite{bhattacharjee2005interval,hettich1993semi,spissueLQW08,mitsos2011global}
for classical surveys on SIPs.

It is harder to solve GSIPs than SIPs.
The methods for SIPs are difficult to be extended to solve GSIPs,
when the set $U(x)$ depends on $x$.
To show the difficulty for doing this, suppose $\mc{F}_0$
is an initial approximation for
the feasible set $\mc{F}$ of \reff{primal-GSIP}
and solve the relaxation problem
\be \label{upper-level subproblem}
\left\{
\begin{array}{cl}
	\min \, & f(x) \\
	\st  \, & x \in \mc{F}_0.
\end{array}
\right.
\ee
Suppose $\hat{x}$ is an optimizer for \eqref{upper-level subproblem}.
Typically, the point $\hat{x}$ is not a feasible point in $\mc{F}$
(otherwise, $\hat{x}$ is a minimizer for the GSIP),
i.e., there exists a point $\hat{u} \in U(\hat{x})$ such that
$g(\hat{x}, \hat{u}) \not \ge 0 .$
Such $\hat{u}$ can be obtained by solving the
optimization problem about $u$:
\be \label{lower-level subproblem}
\left\{
\begin{array}{cl}
	\min \, &  \min\limits_{1 \le j \le s} g_j(\hat{x},u) \\
	\st  \, &  u \in U(\hat{x}).
\end{array}
\right.
\ee
The classical exchange technique (see \cite{hettich1993semi})
is to solve the strengthened relaxation problem
\be \label{upper-level:strength}
\left\{
\begin{array}{cl}
	\min \, & f(x) \\
	\st  \, & g(x, \hat{u}) \ge 0, \,  x \in \mc{F}_0.
\end{array}
\right.
\ee
However, the new optimization~\reff{upper-level:strength}
may not contain any minimizer $x^*$ of \reff{primal-GSIP}
as a feasible point. This is because
$g(x^*, u) \ge 0$ is only required to hold for all $u \in U(x^*)$.
But we may not have $g(x^*, \hat{u}) \ge 0$,
since it is possible that $\hat{u} \not\in U(x^*)$.
This phenomenon is shown in Example~\ref{Alexander 2}.
The classical exchange trick may not solve GSIPs correctly.
Other types of computational methods for SIPs
experience similar difficulties for solving GSIPs.
We refer to \cite{vazquez2008generalized,still1999generalized}
for the classical work.

Some methods transform GSIPs into other equivalent problems
(see \cite{lemonidis2008global} for related work).
When the objective $f$ is convex, $g(x,u)$
is linear in $x$ and is quadratic in $u$,
a branch-and-bound method is proposed in \cite{levitin1998branch} ,
which transforms a GSIP into a finite number of convex programs.
Under certain assumptions, GSIPs can be transformed
to SIPs \cite{still1999generalized}
or some bilevel programs \cite{stein2002GSIPBilevel}.
The GSIP can be converted to min-max optimization \cite{tsoukalas2009global}.
An approach transforming GSIPs into
penalty problems with finitely many constraints
is given in \cite{geletu2004conceptual}.
In computational practice, these equivalent problems
or reformulations are also hard to solve.
It is still a big challenge to solve GSIPs.

%\subsection*{Contributions}
\bigskip \noindent
{\bf Contributions} \,
To solve the GSIP~\eqref{primal-GSIP} efficiently and reliably,
we propose a hierarchy of polynomial optimization relaxations $(P_k)$.
They have the same objective $f(x)$ but their feasible sets
$\mathcal{F}_k$ satisfy the nesting containment
\[
\mathcal{F}_0  \supsetneqq  \mathcal{F}_1  \supsetneqq  \cdots
\mathcal{F}_k  \supsetneqq\cdots \supseteq \mathcal{F}.
\]
The initial relaxation $\mathcal{F}_0$ can be any set containing $\mc{F}$.
At the $k$th stage, we solve the relaxation problem
\[
(P_{k}): \,
\left\{\begin{array}{ll}
	\min & f(x) \\
	\st  & x\in \mathcal{F}_{k}.
\end{array}\right.
\]
Suppose $\hat{x}_k$ is a minimizer of $(P_k)$.
If $\hat{x}_k$ is feasible for \eqref{primal-GSIP},
i.e., $g(\hat{x}_k,u)\geq 0$ for all $u\in U( \hat{x}_k )$,
then it is an optimizer for the GSIP.
If it is not feasible for \eqref{primal-GSIP},
we need to update $(P_k)$ to get
the next tighter relaxation $(P_{k+1})$.

When $\hat{x}_k$ is not feasible for \eqref{primal-GSIP},
there must exist $\hat{u}_k \in U(\hat{x}_k)$ such that
$g_j(\hat{x}_k, \hat{u}_k) < 0$ for some $j$, i.e.,
the inequality $g(\hat{x}_k, \hat{u}_k) \ge 0$ is violated.
When $U(x)$ depends on $x$, we cannot simply select
the new relaxation set $\mc{F}_{k+1}$ as
\[
\mathcal{F}_{k+1}  \,  =  \, \big \{
x\in\mathcal{F}_k : g(x,\hat{u}_k) \ge 0
\big \}.
\]
This is because we may not have $\hat{u}_k \in U(x^*)$
for any minimizer $x^*$ of \reff{primal-GSIP}.
The classical exchange trick does not work for GSIPs.
To fix this issue, we look for a polynomial vector function $q_k(x)$ such that
\[
q_k( \hat{x}_{k} ) \, = \,  \hat{u}_k, \quad
q_k(x) \in U(x) \quad \forall \, x \in X .
\]
Such $q_k$ is called a polynomial extension at $\hat{u}_k$.
It is important to observe that if such $q_k$ exists, then
every feasible point $x$ of \reff{primal-GSIP} must satisfy the inequality
$
g(x, q_k(x) )  \ge  0 .
$
Therefore, if $x^*$ is a minimizer of \reff{primal-GSIP},
then we must have
$
g(x^*, q_k(x^*) )  \ge  0 .
$
Therefore, we can select the new polynomial optimization relaxation $(P_{k+1})$ as
\be \label{contribution feasible set}
\left\{\begin{array}{ll}
\min & f(x)\\
\st  &  g(x, q_k(x) ) \ge 0,  \,
x \in \mathcal{F}_{k} .
\end{array}\right.
\ee
Equivalently, the new feasible set $\mathcal{F}_{k+1}$ is updated as
\[
\mathcal{F}_{k+1} \, \coloneqq \,
\big \{
x\in\mathcal{F}_k : \,  g(x,q_k(x))\geq 0
\big \}.
\]
Continuing to do the above,  we can get a hierarchy of
polynomial optimization relaxations.
This hierarchy gives a method for solving the GSIP.
It produces a sequence of points $\hat{x}_k$,
which serve as approximate optimizers for \eqref{primal-GSIP}.
It is expected that either $\hat{x}_k$ is a minimizer of
\eqref{primal-GSIP} for some $k$
(this case is called the finite convergence),
or an accumulation point of
the sequence $\{ \hat{x}_k \}^{\infty}_{k=0}$
is a minimizer of \eqref{primal-GSIP}
(this case is called the asymptotic convergence).
Under certain assumptions, we can prove that
the hierarchy of polynomial optimization relaxations
has either finite convergence or asymptotic convergence.

This paper is organized as follows.
In Section~\ref{sec:Pre}, we review some preliminaries about polynomial optimization.
In Section~\ref{sec:GSIP}, we show how to use polynomial extensions
to construct the hierarchy of polynomial optimization relaxations.
Its convergence is also shown there.
In Section~\ref{sec:Poly}, we show how to get polynomial extensions.
In Section~\ref{sc:convex}, we discuss how to solve GSIPs
that have convex infinite constraints.
In Section~\ref{sec:Numerical}, various numerical examples for both GSIPs and SIPs
are presented to show the efficiency of the proposed method.
Some conclusions and discussions are made in Section~\ref{sec:COnclu}.

\section{Preliminaries}
\label{sec:Pre}

\noindent
{\bf Notation} \,
Let $\N$ (resp., $\re$) denote the set of nonnegative integers
(resp., real numbers). For $t \in \re$, $\left\lceil t \right\rceil$
denotes the smallest integer greater than or equal to $t$.
For a finite set $A$, we use $\vert A \vert$ to denote its cardinality.
The symbol $\mathbb{R}[x]=\mathbb{R}[x_1,\cdots,x_n]$
denotes the ring of polynomials in
$x \coloneqq (x_1,\cdots,x_n)$ with real coefficients.
For a symmetric matrix $W$, $W\succeq 0$ (resp., $W \succ 0$)
means that $W$ is positive semidefinite (resp., positive definite).
For a vector $u$,  $\|u\|$  denotes the standard Euclidean norm.
The notation $I_n$ denotes $n$-dimensional identity matrix,
and $e$ denotes the vector of all ones
(its length is determined by the operation on it).
For an integer $s>0$, denote the set $[s] \coloneqq \{1,\ldots, s\}$.
For a power $\alpha\coloneqq(\alpha_1,\cdots,\alpha_n)\in\mathbb{N}^n$,
we denote the monomial $x^\alpha\coloneqq x^{\alpha_1}_1\cdots x^{\alpha_n}_n$.
The notation
\[
\mathbb{N}^n_d = \{\alpha\in\mathbb{N}^n :
\alpha_1+\cdots+\alpha_n\leq d\}
\]
denotes the set of monomial powers with degrees at most $d$.
So, $\mathbb{R}^{\mathbb{N}^n_{2d}}$ denotes the space of
real vectors labeled by $\alpha \in\mathbb{N}^n_{2d}$.
The notation $\nabla_x f$ stands for the gradient of $f$ in $x$
and $\nabla_u g$ stands for the gradient of $g$ in $u$.

A polynomial tuple $h=\left(h_{1}, \ldots, h_{m}\right)$
in $\mathbb{R}[x]$ generates the ideal
\[
\text{Ideal}(h) \, \coloneqq \,
h_{1} \cdot \mathbb{R}[x] +\cdots+ h_{m} \cdot \mathbb{R}[x],
\]
the smallest ideal containing all $h_{i}$.
The $2k$-th degree truncation of $\text{Ideal}(h)$ is
\[
\text{Ideal}(h)_{2k}  \coloneqq  h_{1} \cdot \mathbb{R}[x]_{2k-\operatorname{deg}\left(h_{1}\right)}+\cdots+h_{m} \cdot \mathbb{R}[x]_{2k-\operatorname{deg}\left(h_{m}\right)}.
\]
A polynomial $\sigma \in \mathbb{R}[x]$ is said to be a sum of squares (SOS) polynomial if $\sigma=\sigma_{1}^{2}+\cdots+\sigma_{k}^{2}$ for some
$\sigma_{1}, \ldots, \sigma_{k} \in \mathbb{R}[x]$.
The symbol $\Sigma[x]$ denotes the set of all SOS polynomials in $x$.
We refer to \cite{HN08,NiePMI11,NieSOS12} for related work
about SOS polynomials and matrices.
The quadratic module for a polynomial tuple $q = (q_{1}, \ldots, q_{t})$
is ($q_0 \coloneqq 1$)
\[
Q M[q] \coloneqq \left\{ \sum_{i=0}^{t} \sigma_{i} q_{i}\,
\middle\vert\, \sigma_{i} \in \Sigma[x]\right\}.
\]
For a positive integer $k > 0$,
the $2k$th degree truncation of $\text{QM}(q)$ is
\[
Q M[q]_{2k}\coloneqq\left\{\sum_{i=0}^{t}
\sigma_{i} q_{i}\,\middle\vert\, \sigma_{i} \in \Sigma[x],
\operatorname{deg}\left(\sigma_{i} q_{i}\right) \leq 2k\right\}.
\]
The sum $\text{Ideal}(p)+$ $\text{QM}(q)$ is said to be archimedean
if there exists $b \in \operatorname{Ideal}(p)+\operatorname{QM}(q)$
such that the inequality $b(x) \ge 0$ gives a compact set.

\medskip
A vector labelled as
$y    \coloneqq    (y_{\alpha})_{\alpha \in \mathbb{N}^n_{2d}}$
is called a truncated multi-sequences (tms) of degree $2d$.
For $y \in \re^{\mathbb{N}_{2d}^{n}}$, define the operation
\[
\langle p, y \rangle \,  \coloneqq  \,
\sum_{ \alpha \in \mathbb{N}_{2d}^{n} } p_\af y_\af
\qquad \text{for} \quad
p  = \sum_{ \alpha \in \mathbb{N}_{2d}^{n} } p_\af x^\af  .
\]
For a given $p \in\mathbb{R}[x]$ with degree $\deg(p) \le 2d$, let
$t   \coloneqq   \lceil d -\deg(p)/2\rceil$.
Let $L^{(d)}_{p}[y]$ denote the symmetric matrix such that
\[
\langle p \phi^2, y \rangle \quad = \quad
\mathrm{vec}(\phi)^T \cdot  L^{(d)}_{p}[y] \cdot \mathrm{vec}(\phi)
\qquad \text{for all} \quad \phi \in \re[x]_t.
\]
where $\mathrm{vec}(\phi)$ denotes the coefficient vector of $\phi$,
ordered by the graded alphabetical ordering.
In particular, for the case $p=1$, we get the moment matrix
$M_d[y] = L^{(d)}_{1}[y]$, which is equivalent to the expression
\[
M_d[y]  \coloneqq (y_{\alpha+\beta})_{ \alpha, \beta\in\mathbb{N}^n_d } .
\]
For the polynomial $p$, let $\mc{V}_{p}^{(2d)}[y]$ denote the vector such that
\[
\langle p \phi, y \rangle  \, = \,
(\mc{V}_{p}^{(2d)}[y])^T   \mathrm{vec}(\phi).
\]
The $\mc{V}_{p}^{(2d)}[y]$ is called the
localizing vector of $p$, generated by $y$.

\subsection{Moment-SOS relaxations}
\label{ssc:pop}

For convenience of description,
we consider the standard polynomial optimization problem
\be   \label{relaxation}
\left\{
\begin{array}{ll}
	\min\limits_{x \in \re^n} & f(x) \\
	\st    & p_i(x) = 0 \, (i=1,\ldots, m_1) \\
	& p_i(x) \ge 0 \, (i=m_1+1,\ldots, m_2).
\end{array}
\right.
\ee
The Lasserre type Moment-SOS hierarchy
can be applied to solve \reff{relaxation}.
The SOS relaxation is
\be \label{formal relaxation}
\left\{
\begin{array}{cl}
\max \, & \gamma  \\
\st  \, & f-\gamma \in \mathrm{Ideal}\left[p_1,\ldots,p_{m_1}\right]_{2 k}+
\mathrm{QM}\left[p_{m_1+1},\ldots, p_{m_2}\right]_{2 k}.
\end{array}
\right.
\ee
Its dual problem is the moment relaxation
\begin{equation}\label{prelimilary dual sos}
\left\{\begin{array}{cl}
\min \, & \langle f, y\rangle \\
\st  \, &\mathcal{V}_{p_{i}}^{(2 k)}[y]=0\,(i = 1, \ldots, m_1), \\
& L_{p_{i}}^{(k)}[y] \succeq 0\,(i = m_1+1, \ldots, m_2), \\
& M_{k}[y] \succeq 0 ,\\
& y_{0}=1, y \in \mathbb{R}^{\mathbb{N}_{2 k}^{n}}.
\end{array}\right.
\end{equation}
Suppose $y^*$ is an optimizer of $\eqref{prelimilary dual sos}$.
Denote the degree
\[
d_{0}  \coloneqq  \max_{ i \in [m_2] }
\{  \lceil \deg(f) / 2\rceil, \lceil \deg(p_i) / 2\rceil \}.
\]
If there exists a degree $t \in [d_0, k]$ such that
\be \label{flat truncation}
\text{rank} M_{t-d_0}[y^*] = \text{rank} M_t[y^*],
\ee
then we can extract  $r\coloneqq\text{rank} M_t[y^*]$ minimizers.
Generically, flat truncation is a sufficient and necessary condition for
extracting optimizers from Moment-SOS relaxations \cite{nie2013certifying}.
We refer to \cite{HKL20,las01,Lau09,MPO}
for general polynomial optimization theory and methods.

\section{Polynomial optimization relaxations}
\label{sec:GSIP}

This section gives a hierarchy of polynomial optimization relaxations
for solving GSIPs.
It requires to use polynomial extensions.

\subsection{Polynomial extensions}
\label{ssc:PE}

To solve the GSIP \eqref{primal-GSIP} efficiently and reliably,
we relax its infinity constraints
into finitely many ones and then
subsequently update them to get better relaxation.
To achieve this, we need to construct a sequence of
polynomial optimization relaxations $(P_k)$,
which have the same objective $f(x)$ and have the feasible sets
$\mathcal{F}_k$ satisfying the nesting containment
\be \label{nesting:Fk}
\mathcal{F}_0\supsetneqq \mathcal{F}_1\supsetneqq \cdots
\mathcal{F}_k\supsetneqq\cdots \supseteq \mathcal{F}.
\ee
The initial $\mc{F}_0$ can be any set containing $\mc{F}$
(e.g,  $\mathcal{F}_0 = X$).
For each $k=0,1,\ldots$, the GSIP \eqref{primal-GSIP} is relaxed to
\be  \label{Pk equation}
(P_{k})\quad\left\{
\begin{array}{rcl}
	f_k \coloneqq  & \min & f(x)\\
	& \st & x\in \mathcal{F}_k.
\end{array}
\right.
\ee
The minimal value for $(P_k)$ is denoted as $f_k$.
Clearly, the nesting relation \reff{nesting:Fk} implies
\be \label{lowerbd:mono}
f_0 \leq f_1 \leq \cdots f_k \leq  \cdots \leq f^*,
\ee
where $f^*$ denotes the optimal value of \reff{primal-GSIP}.
Suppose $\hat{x}_k$ is a minimizer for \reff{Pk equation}.
To check if $\hat{x}_k$ is feasible for \eqref{primal-GSIP} or not,
we need to solve the optimization problems (for $j=1,\ldots,s$),
\be \label{Qk equation}
(Q_{k,j}) \quad \left\{\begin{array}{rl}
	\hat{g}_{k,j}  \coloneqq
	\min &  g_j(\hat{x}_k,u)  \\
	\st &  u  \in  U(\hat{x}_k).
\end{array}\right.
\ee
The minimal value of  $(Q_{k,j})$ is denoted as $\hat{g}_{k,j}$.
Note that $\hat{x}_k$ is feasible for \eqref{primal-GSIP} if and only if
$g(\hat{x}_k, u) \geq 0$ for all $u  \in  U(\hat{x}_k)$. Therefore,
\[
\hat{x}_k \,\text{is feasible for \eqref{primal-GSIP}}
\quad \Leftrightarrow \quad  \hat{g}_{k,j} \ge  0,
\, j = 1, \ldots, s.
\]
If some $\hat{g}_{k,j} < 0$,  then $\hat{x}_k$
is not feasible for $\eqref{primal-GSIP}$.
Denote the label set
\be \label{labelset:Nk}
\mc{N}_k \, \coloneqq \,
\big\{ j \in [s]: \hat{g}_{k,j} < 0  \big \}.
\ee
Suppose $\hat{u}_{k,j}$ is a minimizer of $(Q_{k,j})$.
Using the classical exchange trick, one may select the new set
$\mathcal{F}_{k+1}$ as
\be \label{feasible set pk+1 original}
\mathcal{F}_{k+1} \, \coloneqq  \,
\{x\in\mathcal{F}_k : \,  g(x,\hat{u}_{k,j}) \ge  0,
j \in \mc{N}_k  \}.
\ee
However, the above choice may result in that $x^* \not\in \mc{F}_{k+1}$
for every optimizer $x^*$ of \eqref{primal-GSIP}.
The reason for this is that the constraining set
$U(x)$ depends on $x$. This can be shown in the following example.

\begin{eg}\label{Alexander 2}
Consider the following GSIP from \cite{articleAngelos2015}:
\[
\left\{
\begin{array}{cl}
	\min \, & -x_1\\
	\st  \, & g \coloneqq u^5-3x^2_2\geq 0\quad\forall u\in U(x),\\
	& x\in X,
\end{array}
\right.
\]
where the constraining sets are
\[
\begin{array}{ll}
	X=[0,1]^2, \quad
	U(x)=\left\{u\in \re^1 \,\middle \vert\,\begin{array}{ll}
		-u\geq 0, u+2\geq 0,\\ u^5+4x^2_1+x^2_2-1\geq 0
	\end{array}\right\}.
\end{array}
\]
The minimizers for $(P_0)$ and $(Q_0)$ are
$\hat{x}_0 \approx  \begin{pmatrix}1.0000,0.4350\end{pmatrix}$
and $\hat{u}_0 \approx -1.2611$ respectively.
The true minimizer for this GSIP is
$x^*=\begin{pmatrix}
	0.5000, 0.0000
\end{pmatrix}$.
However, $g(x^*,\hat{u}_0) \approx -3.1892<0$. If we let
$\mathcal{F}_1 = \{x\in \mathcal{F}_0 :  g(x,\hat{u}_0)\geq 0\}$,
then $x^*\not\in\mathcal{F}_1$ and $\mathcal{F}_1\not\supseteq\mathcal{F}$.
We cannot get the minimizer if $\mc{F}_1$ is selected
as in \reff{feasible set pk+1 original}.
More details for solving this GSIP is referred to
Example~\ref{Box change variable}.
\end{eg}

To fix the above issue, we need to introduce the concept
of polynomial extension (also see \cite{nie2021lagrange}).

\begin{defi} \label{PolyExtenDefi}
For a pair $(\hat{x},\hat{u})$ with $\hat{x} \in X$ and $\hat{u} \in U(\hat{x})$,
a polynomial vector function $q(x)$
is called a polynomial extension of $\hat{u}$ for $U(x)$ if
\be \label{defi:polyexten}
q(\hat{x})=\hat{u}, \qquad
q(x) \in U(x) \quad \forall x \in  X .
\ee
\end{defi}

Recall that $\hat{x}_k$ is a minimizer of \reff{Pk equation}.
Consider the case that $\hat{g}_{k,j} < 0$ and
$\hat{u}_{k,j}$ is a minimizer of $Q_{k,j}$.
Assume $q_{k,j}(x)$ is a polynomial extension of
$\hat{u}_{k,j}$ for $U(x)$, i.e.,
\be \label{polyext:qkj}
q_{k,j}(\hat{x}_k) = \hat{u}_{k,j}, \qquad
q_{k,j}(x)\in U(x) \quad\forall x\in X .
\ee
Observe that $g(x,q_{k,j}(x)) \ge  0$ for every feasible point $x \in \mc{F}$.
Therefore, the new tighter relaxation $\mathcal{F}_{k+1}$
can be selected as
\be \label{update feasible set}
\mathcal{F}_{k+1}  \, \coloneqq  \,
\{x\in\mathcal{F}_k : \,  g(x,q_{k,j}(x)) \ge 0,  \,
j \in \mc{N}_k  \}.
\ee
Continuing to do this, we can get a hierarchy of relaxations $\mc{F}_k$
satisfying the nesting relation \reff{nesting:Fk}.

\subsection{The algorithm}
\label{ssc:alg}

Using polynomial extensions in the above,
we can get the following algorithm for solving GSIPs.

\begin{alg}  \label{algorithm}
For the GSIP \eqref{primal-GSIP}, do the following:
\begin{description}

\item [Step~0] Let $k \coloneqq 0$ and $\mathcal{F}_0 \coloneqq X$
(or any set $\mathcal{F}_0 \supseteq \mc{F}$).

\item [Step~1]
Solve the polynomial optimization relaxation~$(P_k)$ as in \eqref{Pk equation}.
If $(P_k)$ is infeasible, output that \eqref{primal-GSIP}
is infeasible and there are no minimizers.
If it is feasible and has minimizers, compute a minimizer $\hat{x}_k$.

\item [Step~2]
Check whether or not $\hat{x}_k$ is feasible for \eqref{primal-GSIP}.
Solve the optimization problem $(Q_{k,j})$ as in
\eqref{Qk equation}, for each $j=1,\ldots, s$.
If $U(\hat{x}_k) = \emptyset$,
then the constraint that $g(\hat{x}_k, u) \ge 0$
for all $u \in U(\hat{x}_k)$ is automatically satisfied.
For this case, let $\hat{g}_k \coloneqq + \infty$.
If $U(\hat{x}_k) \ne \emptyset$,
compute its optimal value $\hat{g}_{k,j}$
and a minimizer $\hat{u}_{k,j}$ (if any) for $(Q_{k,j})$.
Then let
$
\hat{g}_k   \coloneqq
\min\limits_{ j \in [s] } \, \hat{g}_{k,j} .
$

\item [Step~3]
If $\hat{g}_k \ge 0$, output that $\hat{x}_k$ is a minimizer
for \eqref{primal-GSIP} and stop.
If $\hat{g}_k <  0$, construct a polynomial extension
$q_{k,j}(x)$ satisfying \reff{polyext:qkj}
for each $j \in \mc{N}_k$,
where $\mc{N}_k$ is the label set as in \reff{labelset:Nk}.

\item [Step~4]
Update the set $\mathcal{F}_{k+1}$ as in \eqref{update feasible set}.
Let $k \coloneqq k+1$ and go to Step~1.

\end{description}

\end{alg}

We would like to remark that Algorithm~\ref{algorithm}
is in the framework of the {\it conceptual} exchange method
of Still's work \cite[Algorithm~2]{still2001generalizedaspects}.
General nonlinear functions are considered in \cite{still2001generalizedaspects}.
There are no efficient methods/packages for computing
global minimizers for general nonlinear optimization.
Moreover, computing continuous functions satisfying
feasible sets given by infinitely many constraints is also difficult.
The conceptual exchange method is
generally hard to be implemented in computational practice.
However, when GSIPs are given by polynomials,
the optimization problems $(P)_k$ and $(Q_{k,j})$
can be solved globally (by Moment-SOS relaxations, see section~\ref{ssc:pop}),
and the polynomial extensions
can be computed for many frequently appearing feasible sets
(see Section~\ref{sec:Poly}).
This is our major motivation for studying GSIPs of polynomials.

To check if $\hat{x}_k$ is feasible  for \eqref{primal-GSIP} or not,
the optimization problem $(Q_{k,j})$ in \eqref{Qk equation}
need to be solved very accurately.
We can apply the tight relaxation method in
\cite{nie2019tight} to solve them.
In Step~2, if $(Q_{k,j})$ is unbounded below, i.e., $\hat{g}_{k,j} = - \infty$,
we can just choose any point $\hat{u}_{k,j}$
such that $g_j(\hat{x}_k, \hat{u}_{k,j}) < 0$.
In computational practice, we often select a priori accuracy parameter
$\eps > 0$ (e.g, $\eps = 10^{-6}$), to detect if
$\hat{x}_k$ is feasible for \eqref{primal-GSIP} or not.
If $\hat{g}_k  \ge  -\epsilon$,
we can view $\hat{x}_k$ as a feasible point,
i.e., it is an accurate minimizer for the GSIP.

The polynomial optimization relaxation $(P_k)$ can be solved
by the Lasserre type Moment-SOS relaxations
(see subsection \ref{ssc:pop}).
Under the archimedean condition,
it has the asymptotic convergence (see \cite{las01}).
Moreover, under some classical optimality conditions,
the finite convergence occurs (see \cite{nie2013exact}).
The convergence can be detected by the flat truncation condition
(see \cite{nie2013certifying}).

\subsection{Convergence properties}
\label{ssc:cvg}

We now study the convergence of Algorithm~\ref{algorithm}.
First, if Algorithm~\ref{algorithm} terminates at a certain loop,
the computed solution is a minimizer
for the GSIP\eqref{primal-GSIP}.

\begin{thm}[Finite convergence]
	\label{Convergence_theorem}
	In Algorithm~\ref{algorithm}, if $\hat{g}_k \ge 0$ at the $k$th loop,
	then the minimizer $\hat{x}_k$ for $(P_k)$
	is also a minimizer for \eqref{primal-GSIP}.
\end{thm}
\begin{proof}
The feasible sets $\mathcal{F}_k$
satisfy the nesting relation \reff{nesting:Fk}
and the optimal values $f_k$ satisfy \reff{lowerbd:mono}.
If $\hat{g}_k \ge 0$ at the $k$th loop,
it means that $g(\hat{x}_k, u)\ge 0$ for all $u\in U(\hat{x}_k)$,
i.e., $\hat{x}_k$ is feasible for \eqref{primal-GSIP}.
So, $f_k = f(\hat{x}_k) \ge f^*$, hence $f_k=f^*$
and $\hat{x}_k$ is a minimizer for \eqref{primal-GSIP}.
\qed
\end{proof}

Second, we study the asymptotic convergence of Algorithm~\ref{algorithm}.
Let $x^*$ be an accumulation point of the sequence
$\{ \hat{x}_k \}_{k=1}^\infty$.
We consider the value function
\be \label{valfun:v(x)}
v(x) \, \coloneqq  \, \min\limits_{ j \in [s] }
\min\limits_{  u\in U(x) } g_j(x,u) .
\ee
People often assume $v(x)$ is continuous at $x^*$.
This can be ensured under certain conditions
(see \cite{clarke1998nonsmooth,clarke1990optimization,guo2014sensitivity,lucet2001sensitivity}).
For convenience of notation, let $\hat{u}_k$
denote the point $\hat{u}_{k,j}$ such that
$g_j(\hat{x}_k, \hat{u}_{k,j} ) = \hat{g}_k$
and its polynomial extension $\hat{q}_{k,j}(x)$ of $\hat{u}_{k,j}$
is also denoted as $\hat{q}_k(x)$.

\begin{thm}[Asymptotic convergence]
\label{Asympotic_theorem}
Let $\{ \hat{x}_k \}_{k=0}^\infty$ be the sequence
produced by Algorithm~\ref{algorithm}. Suppose
$x^*$ is an accumulation point of this sequence.
Assume the following:
\bit

\item [(i)]
In the $k$th loop, the problems $(P_k)$ and $(Q_{k,j})$
have minimizers $\hat{x}_k$ and $\hat{u}_{k,j}$ respectively;

\item [(ii)]
Suppose $(x^*,u^*)$ is an accumulation point of
$\{(\hat{x}_k,\hat{u}_k)\}^{\infty}_{k=0}$ and
the value function $v(x)$ as in \reff{valfun:v(x)} is continuous at $x^*$;

\item [(iii)]
The polynomial extensions $q_{k,j}(x)$ are equicontinuous at $x^*.$\footnote{
This means that for every $\eps>0$, there exists $\dt >0$
such that $\| q_{k,j}(x) - q_{k,j}(x^*) \| < \eps$
for all $\| x - x^* \| < \dt$ and for all $k,j$. }

\eit
Then, we have $v(x^*) = 0$ and $x^*$
is a minimizer for the GSIP \eqref{primal-GSIP}.
\end{thm}
\begin{proof}
Since $x^*$ is an accumulation point of the sequence
$\{ \hat{x}_k \}_{k=1}^\infty$, there exists a subsequence
such that $\lim\limits_{l \rightarrow\infty}\hat{x}_{k_l}= x^*$.
Since $f$ is a polynomial function,
the monotonicity relation \reff{lowerbd:mono} implies
\[
f(x^*) = \lim\limits_{l \rightarrow \infty} f(\hat{x}_{k_l})
= \lim\limits_{k_l\rightarrow \infty} f_{k_l} \leq f^*.
\]
If Algorithm~\ref{algorithm} does not stop at the $k_l$th loop, we have
\[
v(\hat{x}_{k_l}) \, < \, 0 .
\]
Letting $l  \rightarrow  \infty$, we get $v(x^*)\leq 0$,
since $v(x)$ is continuous at $x^*$.

Next, we pick an arbitrary $k^\prime \in \N$.
For every $k_l > k^\prime$,  the feasibility of
the polynomial optimization relaxation $(P_{k_l})$ implies that
\[
g(\hat{x}_{k_l}, \hat{q}_{k^\prime} (\hat{x}_{k_l})) \ge 0 .
\]
Letting $l \rightarrow \infty$, we then get
$
g(x^*,  \hat{q}_{ k^\prime } (x^*) ) \ge  0 .
$
Since $k^\prime$ can be arbitrary, the above implies that
\be \label{gx*qk>=0}
g(x^*,  \hat{q}_{k}(x^*) ) \ge 0 ,
\quad \mbox{for all} \quad k \in \N .
\ee
Note that $\hat{u}_{k_l}$ denotes the point $\hat{u}_{k_l,j}$ such that
$g_j(\hat{x}_{k_l}, \hat{u}_{k_l,j} ) = \hat{g}_{k_l}$, so
\[
g_j(\hat{x}_{k_l}, \hat{u}_{k_l} ) \, = \, g_j(\hat{x}_{k_l}, \hat{u}_{k_l, j} ) .
\]
Since $\hat{u}_{k_l,j}$ is the minimizer of
$g_j(\hat{x}_{k_l}, u)$ over $u \in U(\hat{x}_{k_l})$,
we have
\be  \label{vhat_xkl}
v( \hat{x}_{k_l}) = \min\limits_{j \in [s]} g_j(\hat{x}_{k_l}, \hat{u}_{k_l, j} )
= \min\limits_{j \in [s]} g_j(\hat{x}_{k_l}, \hat{u}_{k_l} ).
\ee
Therefore, \reff{gx*qk>=0} and \reff{vhat_xkl} imply that
\begin{eqnarray*}
	v(x^*) & = & v( \hat{x}_{k_l} ) + \big( v(x^*) -  v( \hat{x}_{k_l} )  \big) \\
	& \ge & \Big( \min\limits_{j \in [s]} g_j( \hat{x}_{k_l}, \hat{u}_{k_l} ) -
	\min\limits_{j \in [s]} g_j(x^*, \hat{q}_{k_l} (x^*) ) \Big)
	+ \big( v(x^*) -  v( \hat{x}_{k_l} ) \big).
\end{eqnarray*}
By the polynomial extension and equicontinuity, we have
\[
u^* = \lim\limits_{l \to \infty} \hat{u}_{k_l} =
\lim\limits_{l \to \infty} \hat{q}_{k_l} (\hat{x}_{k_l}) =
\lim\limits_{l \to \infty} \hat{q}_{k_l} (x^*),
\]
\[
\lim\limits_{l \to \infty}
\Big(  \min\limits_{j \in [s]} g_j( \hat{x}_{k_l}, \hat{u}_{k_l} ) -
\min\limits_{j \in [s]} g_j(x^*, \hat{q}_{k_l} (x^*) \Big)
\,  = \,  0.
\]
Since $v(x)$ is continuous at $x^*$,
the above implies $v(x^*) \ge 0$.

Combining the above, we have $v(x^*)=0$,
so $g(x^*,u)\geq 0$ holds for all $u\in U(x^*)$.
This implies that $x^*$ is feasible for the GSIP \eqref{primal-GSIP},
i.e., $x^*\in\mathcal{F}$. Hence,  $f(x^*) \geq f^*$
and we must have $f(x^*)=f^*$. In other words, $x^*$
is a minimizer for the GSIP \eqref{primal-GSIP}.
\qed
\end{proof}

When the sequence $\{ \hat{x}_k \}^{\infty}_{k=0}$ is bounded,
it always has an accumulation point.
This is the case if either $X$ is compact or the objective
$f$ is coercive, i.e., the sublevel set $\{x \in X: f(x)\leq \upsilon\}$
is bounded for each value $\upsilon$.
The equicontinuity condition in Algorithm~\ref{Asympotic_theorem} (iii)
can be ensured if the polynomial extension functions $q_{k,j}(x)$
have uniformly bounded degrees and bounded coefficients.
In Section~\ref{sec:Poly}, we will see that this is the case
for many common feasible sets.

\section{Polynomial Extensions}
\label{sec:Poly}

This section discusses how to compute polynomial extensions.
Suppose $\hat{x}$ is an approximate optimizer obtained from
solving a polynomial optimization relaxation.
We consider the case that $\hat{x}$ is not feasible
for the GSIP \reff{primal-GSIP},
i.e., there exists $\hat{u} \in U(\hat{x})$
such that $g_j(\hat{x}, \hat{u}) < 0$ for some $j\in [s]$.
For the special case of SIP,
the polynomial function $q(x)$ can be selected to be
the constant $\hat{u}$. In the following,
we study how to compute a polynomial function
$q(x)$ satisfying \reff{defi:polyexten} in Definition~\ref{PolyExtenDefi}
for more general cases.

First, we look at some common feasible sets for which polynomial extensions
can be given explicitly. They are discussed in \cite{nie2021lagrange}.

\bigskip \noindent
{\bf Box type constraints} \,
Suppose the infinity constraining set is
\[	
U(x) \, =  \, \left\{u \in \mathbb{R}^{p} : \,
l(x) \leq u \leq w(x) \right\},
\]
where $l(x) \, \coloneqq \, (l_{1}(x), \cdots, l_{p}(x))$ and
$w(x) \, \coloneqq \,  (w_{1}(x), \cdots, w_{p}(x))$
are given polynomial vectors.
Suppose $l(x) \le w(x)$ on $X$.
For the given pair $(\hat{x}, \hat{u})$,
with $\hat{u} = (\hat{u}_1, \ldots, \hat{u}_p)$,
we can choose the polynomial extension
$q= (q_1,\ldots, q_p)$ as
\be \label{PE:box}
q_j(x) = \sigma_{j} l_{j}(x)+\left(1-\sigma_{j}\right) w_{j}(x),
\quad j=1, \ldots, p,
\ee
where each scalar
\[
\sigma_{j} \, \coloneqq \,  (w_{j}(\hat{x})-\hat{u}_j)
/ (w_{j}(\hat{x})-l_{j}(\hat{x})) \in[0,1].
\]
For the case $w_j(\hat{x})-l_{j}(\hat{x}) = 0$,
we just let $\sigma_{j}=0$.
Then, for each $j = 1, \ldots, p$,
\[
q_j(\hat{x}) = \sigma_{j} l_{j}(\hat{x})+
\left(1-\sigma_{j}\right) w_{j}(\hat{x})=\hat{u}_j .
\]
Clearly, $q(x) \in U(x)$ for all $x\in X$.
Such $q(x)$ satisfies \reff{defi:polyexten}.

\bigskip \noindent
{\bf Simplex type constraints} \
Suppose the infinite constraining set is
\[
U(x) \, = \, \left\{u \in \mathbb{R}^{p} :\,  l(x) \leq u, \,
e^{T} u \leq w(x) \right\},
\]
where $l(x) \, \coloneqq \, (l_{1}(x), \cdots, l_{p}(x))$
is a polynomial vector and
$w(x)$ is a scalar polynomial.
Here, $e$ is the vector of all ones.
Suppose $w(x) - e^T l(x) \ge 0$ on $X$.
For given $(\hat{x}, \hat{u})$,
with $\hat{u} = (\hat{u}_1, \ldots, \hat{u}_p)$,
we can choose the polynomial extension
$q = (q_1, \ldots, q_p )$ as
\be \label{PE:simplex}
q_j(x)  = \mu_{j} \left(w(x)-e^{T} l(x)\right)+l_{j}(x), \quad j=1, \ldots, p,
\ee
where each scalar
\[
\mu_j = (\hat{u}_j-l_j(\hat{x}))/(w(\hat{x})-e^{T} l(\hat{x})) \ge 0.
\]
For the case $w(\hat{x})-e^{T} l(\hat{x}) = 0$,
one can just let $\mu_{j}=0$.
Note that $\sum_{i=}^{p} \mu_i \leq 1$ and
\[
q_j(\hat{x}) = \mu_{j} \cdot ( w(\hat{x})-e^{T} l(\hat{x}) )
+l_{j}(\hat{x})=\hat{u}_j .
\]
Since $e^T l(x) \le w(x)$ on $X$, we have $q(x)\geq l(x)$ and
\[
e^{T} q(x)=e^{T} l(x) (1-\sum_{j=1}^{p} \mu_{j} ) +
(\sum_{j=1}^{p} \mu_{j} ) w(x) \leq w(x).
\]
Therefore, $q(x) \in U(x)$ for all $x\in X$.	
Such $q(x)$ satisfies \reff{defi:polyexten}.

\bigskip \noindent
{\bf Ball type constraints} \
Suppose the infinite constraining set is
\[
U(x)  =  \left\{u \in \mathbb{R}^{p} :
l(x) \leq\|u-a(x)\| \leq w(x)\right\},
\]
where $a(x) = (a_{1}(x), \ldots, a_{p}(x))$ is a polynomial vector
and $l(x), w(x)$ are given polynomials such that
$0 \leq l(x) \leq w(x)$ on $X$.
We can choose the polynomial extension $q$ such that
\be \label{PE:ball}
q(x)   \coloneqq   a(x)+\tilde{q}(x) v, \quad
\tilde{q}(x) = \mu_{1} l(x) + \mu_{2} w(x)
\ee
where the scalars $\mu_{1}, \mu_{2}$ and the vector $v$
are selected such that
\[
\baray{c}
\|\hat{u}-a(\hat{x})\| = \mu_{1} l(\hat{x})+\mu_{2} w(\hat{x}), \\
\mu_{1} \ge 0, \,  \mu_{2} \ge 0, \, \mu_{1} + \mu_{2} = 1, \\
v = (\hat{u} - a(\hat{x}))/\|\hat{u}-a(\hat{x})\|.
\earay
\]
For the case $\hat{u}-a(\hat{x}) = 0$,
we just let $v = e/\sqrt{p}$.
Since $\tilde{q}(\hat{x}) = \|\hat{u} - a(\hat{x})\|$
and $\| v \| = 1$, one can check that
\[
\hat{u}-q(\hat{x}) = (\hat{u}-a(\hat{x}))-(q(\hat{x})-a(\hat{x}))=
(\hat{u}-a(\hat{x}))-\tilde{q}(\hat{x}) v  =  0,
\]
\[
\|q(x)-a(x)\|=\left\|\tilde{q}(x) v \right\| =
\left|\tilde{q}(x)\right| \cdot\| v \|=
\left|\tilde{q}(x)\right| .
\]
Since $0 \leq l(x) \leq w(x)$ on $X$, we have
\[
l(x) \leq\|q(x)-a(x)\| \leq w(x).
\]
The above $q(x)$ satisfies \reff{defi:polyexten}.

\bigskip
%%\subsection{ellipsoidal case} \label{Ellipse}
Beyond the above, we consider the case that
the infinity constraining set is given
by a parameterized ellipsoid
\[
U(x) \, =  \, \big \{ u \in \mathbb{R}^p : \,
(u-a(x))^T (D(x)^T D(x))^{-1} (u-a(x)) \leq 1  \big \} .
\]
Here, $a(x)$ is the center vector, and $D(x)\in\mathbb{R}^{p\times p}$
is a polynomial matrix that is nonsingular on $X$.
We can choose the polynomial extension as
\be \label{polyext:Ellipsoid}
q(x) = a(x)+D(x)^T D(\hat{x})^{-T} (\hat{u} - a(\hat{x})).
\ee
One can see that
\[
q(\hat{x}) = a(\hat{x})+D(\hat{x})^T
D(\hat{x})^{-T}(\hat{u}-a(\hat{x}))=\hat{u}.
\]
Moreover, we have $q(x) \in U(x)$ for all $x$, since
\[
\baray{rl}
& (q(x)-a(x))^T (D(x)^T D(x))^{-1} (q-a(x)) \\
= & (\hat{u}-a(\hat{x}))^T D(\hat{x})^{-1} D(\hat{x})^{-T} (\hat{u}-a(\hat{x}))  \\
= & (\hat{u}-a(\hat{x}))^T (D(\hat{x})^{T} D(\hat{x}))^{-1} (\hat{u}-a(\hat{x})) \\
\le & 1.
\earay
\]
The last inequality is due to that $\hat{u} \in U(\hat{x})$.
The $q(x)$ as in \reff{polyext:Ellipsoid}
is a satisfactory polynomial extension.

In the following, we propose a numerical method for
obtaining polynomial extensions.
Recall that the infinity constraining set $U(x)$ is given by
\[
h_i(x, u) \, = \, 0 \,(i \in \mc{I}_2), \quad
h_i(x, u) \, \ge \, 0 \,(i \in \mc{J}_2).
\]
When $h_i(x, u)$ is nonlinear in $u = (u_1, \ldots, u_p)$,
how to obtain polynomial extensions
is generally highly challenging.
This question is mostly open, to the best of the authors's knowledge.
Here, we consider a typical case that
all $h_i(x, u)$ are linear in $u$
(but they can be nonlinear in $x$), i.e.,
\[
h_i(x, u) \, = \, h_{i,0}(x) + (h_{i,1}(x))^T u,
\]
where each $h_{i,1}$ is a $p$-dimensional polynomial vector.
Recall that the constraining set $X$ is given as in \reff{set:X}.
As required in Definition~\ref{PolyExtenDefi}
for given $(\hat{x}, \hat{u} )$,
the polynomial extension vector $q = (\phi_1(x), \ldots, \phi_p(x) )$
can be obtained by solving the following feasibility problem
(for a given degree $k$)
\be \label{polyext:q:num}
\left\{ \baray{l}
(\phi_1(\hat{x}), \ldots, \phi_p(\hat{x}) )  =  \hat{u}, \\
h_i(x, q(x))  \in  \mathrm{Ideal}[c_{eq}]_{ 2k }
\quad \mbox{for} \quad i \in \mc{I}_2,   \\
h_i(x, q(x))  \in
\mathrm{Ideal}[c_{eq}]_{ 2k } +
\mathrm{QM}[c_{in}]_{ 2k }
\quad \mbox{for} \quad i \in \mc{J}_2,   \\
\phi_1, \ldots, \phi_p \in  \re[x]_{\ell} .
\earay
\right.
\ee
In the above, $c_{eq} = (c_i)_{i \in \mc{I}_1}$
is the tuple of equality constraining polynomials for $X$
and $c_{in} = (c_j)_{j \in \mc{J}_1}$
is the one for inequality constraints, and
\[
\ell \, = \, \min\limits_{i \in  \mc{I}_2 \cup  \mc{J}_2  } \, \{ 2k - \deg(h_{i,1}) \}.
\]
The problem \reff{polyext:q:num}
is a linear convex conic system about coefficients of $q$.
The second constraint implies $h_i(x, q(x))=0$
for all $x \in X$ and $i \in \mc{I}_2$,
and the third one implies $h_i(x, q(x)) \ge 0$
for all $x \in X$ and $i \in \mc{J}_2$.
Therefore, if $q$ satisfies \reff{polyext:q:num},
then $q$ is a polynomial extension of $\hat{u}$ at $\hat{x}$.
We would like to remark that \reff{polyext:q:num}
can be solved as a semidefinite program, for the given degree $k$.
Its size can be estimated as follows.
For each of $\phi_1,\ldots, \phi_p$, there are
$\binom{n+\ell}{\ell}$ coefficients.
For each $i \in  \mc{I}_2 \cup  \mc{J}_2$, the membership there
poses $\binom{n+2k}{2k}$ linear equations about coefficients of $q$.
Each polynomial in $\mathrm{QM}[c_{in}]_{ 2k }$
can be represented some psd matrix variables whose lengths
are at most $\binom{n+k}{k}$.
We refer to \cite[Sec.~2.5]{MPO} for more details
about semidefinite programs arising from \reff{polyext:q:num}.
The feasibility system~\reff{polyext:q:num}
can be implemented in \texttt{YALMIP} \cite{lofberg2004yalmip}
and solved by \texttt{SeDuMi} \cite{strum2001sedumi}.

\begin{eg} \label{exmp:4.5}
Consider $X,U(x)$ are given by
\begin{align*}
& \hat{x} =(1/2, 1),\ \hat{u}=(1/3, 2/3),\, c(x) = (x_2-x_1,x_1,2-x_1-x_2)\geq 0,\\
& h(x,u) = (x_2-x_1u_1-x_2u_2,u_1+x_1u_2-x_1x_2,x_1+3u_1-u_2,u_1-x^2_1u_2)\geq 0.
\end{align*}
A polynomial extension $q = (\phi_1, \phi_2)$ satisfying \reff{polyext:q:num} is
\[
\begin{aligned}
& \phi_1  =  0.2705-0.4539x_1-0.0395x_2+0.2690x^2_1-0.4113x_1x_2+0.0564x^2_2 ,\\
& \phi_2  = 0.7102-0.3609x_1+0.1369x_2.
\end{aligned}
\]
One can check that $q(\hat{x})=\hat{u}$ and for every $x \in X$,
\[
h(x,q(x)) \, \ge \,
( 3.9\cdot 10^{-10}  , 1.0\cdot  10^{-8}   ,0.0151  ,  0.0072) \geq 0.
\]
The system~\reff{polyext:q:num} is implemented
in \texttt{YALMIP} \cite{lofberg2004yalmip}
and solved by \texttt{SeDuMi} \cite{strum2001sedumi}.
A numerical certificate for the above $q$
satisfying \reff{polyext:q:num} is given in Appendix~\ref{AppendixC}.
\end{eg}

\section{Convex Infinity Constraints}
\label{sc:convex}

This section discusses GSIPs that have convex infinity constraints.
Such a GSIP can be solved by
a single polynomial optimization relaxation.

Recall that
$
g(x,u) =  \big( g_1(x,u)), \ldots, g_s(x,u) \big)
$
and the infinity set $U(x)$ is given as in \reff{U(x)set}.
Suppose each $g_j(x,u)$ is convex in $u$ for each fixed $x \in X$.
Assume $U(x)$ is a convex set\footnote{
	When $U(x)$ is convex, its equality constraints
	are typically given by linear equations,
	which can be dropped by eliminating variables.
}
such that
\be \label{convexU(x)}
U(x)  \, =  \,  \left\{ u\in\mathbb{R}^p  :
h_1(x,u) \ge 0, \ldots, h_m(x,u) \ge 0
\right\} ,
\ee
where each $h_i(x,u)$ is concave in $u$ for fixed $x \in X$.
For convenience, denote
\[
h(x,u) \, =\, \bbm h_1(x,u) & \cdots & h_m(x,u) \ebm^T.
\]
In this section, the functions $g_j(x,u), h_i(x,u)$ are
not necessarily polynomials.
They can be rational functions in $(x,u)$.
Under certain constraint qualifications, the infinity constraint
can be equivalently expressed through KKT conditions.

Consider the optimization problem in $u$
(for fixed $x \in X$ and $j\in [s]$)
\be \label{min:gjxu:cvx}
\left\{
\begin{array}{rl}
	\min\limits_{u} &  g_j(x,u) \\
	\st &  h(x,u) \ge 0 .
\end{array} \right.
\ee
Under certain constraint qualification conditions,
if $z_j$ is a minimizer of \reff{min:gjxu:cvx},
there exists the Lagrange multiplier vector
$\lmd_j = (\lmd_{j,1}, \ldots, \lmd_{j,m})$ such that
\be \label{cvxKKT:mingjxu}
\boxed{
\begin{array}{c}
\nabla_u h(x, z_j) \lambda_j
= \nabla_u g_j( x, z_j ),  \\
0 \le h( x, z_j ) \perp \lambda_j \ge  0  .
\end{array}
}
\ee
In the above,  $\nabla_u h(x, z_j) =
\bbm  \nabla_u h_1(x, z_j) & \cdots & \nabla_u h_m(x, z_j) \ebm$
is the transposed Jacobian of $h(x,u)$ at $u=z_j$,
and $\perp$ denotes the perpendicularity relation.
Under the convexity assumption, the constraint that
$g_j(x, u) \ge 0$ for all $u \in U(x)$
is equivalent to that
\[
g_j(x, z_j) \ge 0,   \quad z_j
\,\, \text{satisfies} \quad \reff{cvxKKT:mingjxu}.
\]
If the minimum value of \reff{min:gjxu:cvx}
is achieved at a KKT point for every $x \in X$,
the GSIP \reff{primal-GSIP} is then equivalent to
\be \label{cvx:GSIP}
\left\{
\begin{array}{cl}
\min\limits_{ \substack{ x, z_1, \ldots, z_s \\ \lmd_1, \ldots, \lmd_s }  } & f(x)\\
\st  & \nabla_u h(x, z_j) \lambda_j
= \nabla_u g_j( x, z_j ), \, j=1, \ldots, s,  \\
& 0 \le h( x, z_j ) \perp \lambda_j \ge  0,  \, j= 1, \ldots, s, \\
&  g_j(x, z_j) \geq 0, \,  j=1, \ldots, s,  \\
&    x \in X .
\end{array}
\right.
\ee
The optimization problem \reff{cvx:GSIP} has new variables
$z_1,\ldots, z_s, \lmd_1, \ldots, \lmd_s$,
in addition to the original variable $x$.
Therefore, we get the following theorem.

\begin{thm}
Suppose $U(x)$ is given as in \reff{convexU(x)}.
For each $x\in X$, assume every $g_j(x,u)$
is convex in $u$ and every $h_j(x,u)$ is concave in $u$.
Assume the optimal value of \reff{min:gjxu:cvx}
is achieved at a KKT point $z_j$ for each $x \in X$.
Then, the GSIP \reff{primal-GSIP} is equivalent to \reff{cvx:GSIP}.
\end{thm}

Under certain assumptions, Lagrange multipliers $\lmd_j$
can be eliminated. Denote that (for $j=1, \ldots, m$)
\be \label{C(x,u):eta}
H(x,u) = \begin{bmatrix}
\nabla_u h_1  & \nabla_u h_2 &\cdots&  \nabla_u h_m \\
h_1 &0&\cdots&0\\
0&h_2 &\cdots&0\\
\vdots&\vdots&\ddots&\vdots\\
0&0&\cdots&h_m
\end{bmatrix},
\quad
\eta_j(x,u)  =
\begin{bmatrix}
\nabla_u g_j\\
0\\
\vdots\\
0
\end{bmatrix}.
\ee
The KKT system~\reff{cvxKKT:mingjxu} implies that
\[
H(x,u) \lmd_j \, = \, \eta_j(x,u),
\]
If there exists a matrix polynomial $T(x,u)$ such that
\begin{equation}\label{RLME equation}
	T(x,u)H(x,u) \, = \, \psi(x,u) \cdot I_m,
\end{equation}
for a nonzero polynomial $\psi(x,u)$, then
\[
\lmd_j =  \frac{1}{\psi(x,u)} T(x,u) H(x,u) \lmd_j
\, = \, \frac{1}{\psi(x,u)}  T(x,u) \eta_j(x,u).
\]
We remark that for priori degrees for $T(x,u)$ and $\psi(x,u)$,
the equality \reff{RLME equation} gives linear equations about their coefficients.
In other words,  $T(x,u)$ and $\psi(x,u)$ can be determined
by linear equations once their degrees are given.
It is typically quite hard to estimate the degrees of
$T(x,u)$ and $\psi(x,u)$ satisfying \reff{RLME equation}.
However, for special constraints (e.g., ball, simplex, boxes),
we can get $T(x,u)$ and $\psi(x,u)$ explicitly.
We refer to \cite{nie2019tight} and \cite[Sec.~6.2]{MPO}
for more details about this topic.

When $H(x,u)$ is not identically zero, such a matrix polynomial
$T(x,u)$ always exists (see \cite{nie2021rational}).
If we denote the rational vector
\be \label{RLME:pj}
p_j(x,u) = \frac{1}{\psi(x,u)}  T(x,u) \eta_j(x,u),
\ee
the Lagrange multiplier can be expressed as $\lmd_j = p_j(u)$.
Then, the optimization problems~\reff{primal-GSIP}
and \reff{cvx:GSIP} are equivalent to
\be \label{cvx:GSIP:lme}
\left\{
\begin{array}{cl}
\min\limits_{ \substack{ x, z_1, \ldots, z_s  }  } & f(x)\\
\st  & \nabla_u h(x, z_j) p_j(x, z_j)
= \nabla_u g_j( x, z_j ), \, j=1, \ldots, s,  \\
& 0 \le h( x, z_j ) \perp p_j(x,z_j) \ge  0,  \, j=1, \ldots, s,  \\
&  g_j(x, z_j) \geq 0,  \, j=1, \ldots, s,  \\
&   x \in X .
\end{array}
\right.
\ee
The optimization problem~\reff{cvx:GSIP:lme} does not have
Lagrange multiplier variables. It
is about polynomials or rational functions,
so it can be solved by Moment-SOS relaxations.
We refer to \cite{BHL16,henrion2009gloptipoly}
and \cite[Section~5]{nie2021rational}
for how to solve rational optimization problems.
We also have the following theorem.

\begin{thm}
Suppose $U(x)$ is given as in \reff{convexU(x)}.
For each $x\in X$, assume every $g_j(x,u)$
is convex in $u$ and every $h_j(x,u)$ is concave in $u$.
Assume the optimal value of \reff{min:gjxu:cvx}
is achieved at a KKT point $z_j$ for each $x \in X$
and the Lagrange multiplier vector $\lmd_j = p_j(z_j)$
for $p_j$ as in \reff{RLME:pj}.
Then, the GSIP \reff{primal-GSIP} is equivalent to
the optimization problem~\reff{cvx:GSIP:lme}.
\end{thm}

Computational results for GSIPs with convex infinity constraints
are shown in Examples~\ref{LME_example1} and \ref{LME_example2}.

\section{Numerical Experiments}
\label{sec:Numerical}

This section reports numerical examples to show how
Algorithm~\ref{algorithm} solve GSIPs. The computations are implemented
in MATLAB R2022b on a laptop equipped with a 10th Generation Intel® Core™ i7-10510U
processor and 16GB memory. The Moment-SOS relaxations are
implemented by the software  \texttt{Gloptipoly}\cite{henrion2009gloptipoly},
which calls the software \texttt{SeDuMi} \cite{strum2001sedumi}
for solving the corresponding semidefinite programs.

For the GSIP \eqref{primal-GSIP}, we use $x^*$ and $f^*$ to denote the global minimizer
and the global minimum value respectively. The loops in Algorithm~\ref{algorithm}
are labelled by $k$. In the $k$th loop, $\hat{x}_k$ denotes the minimizer of
$(P_k)$ and $f_k$ denotes the minimum value.
Due to numerical errors, we consider $\hat{x}_k$
as an accurate minimizer for \eqref{primal-GSIP}
if all $\hat{g}_{k,j} \ge - \eps$,
for a positive small value $\eps$ (say, $\epsilon=10^{-6}$).
Thus, the inequality $\hat{g}_{k} \ge - \eps$ serves as the termination criterion.
Then point $\hat{u}_k$ denotes the minimizer $\hat{u}_{k,j}$ for the optimization $(Q_{k,j})$
such that $g_j(\hat{x}_k, \hat{u}_{k,j})  = \hat{g}_k$.
When the algorithm terminates at the $k$th loop,
we output the minimizer $x^*\coloneqq\hat{x}_k$
for the GSIP and the minimum value $f^*\coloneqq f_k$.
Consumed computational time is denoted as {\text time}.
All computational results are displayed with $4$ decimal digits.

\subsection{SIP examples}

Since SIPs are special cases of GSIPs,
Algorithm~\ref{algorithm} can be applied to solve SIPs.
This subsection reports such examples.

\begin{eg}\label{SIP_ex_1}
Consider the following min-max SIP from \cite{pang2020discretization}:
\begin{equation*}
\left\{
\begin{array}{cl}
\min\limits_{x\in X} \, & \max\limits_{u\in U}
         5x^2_1+5x^2_2-u^2_1-u^2_2+x_1(-u_1+u_2+5)+x_2(u_1-u_2+3)\\
\st \, &  \begin{bmatrix}
	0.2-x^2_1-u^2_1\\
	0.1-x^2_2-u^2_2
\end{bmatrix}\geq 0\quad \forall u\in U=[-0.2,0.2]^2,\\
& x \in X=[-100,100]^2.\\
\end{array}
\right.
\end{equation*}
This min-max semi-infinite problem is equivalent to an SIP problem:
\begin{equation*}
\left\{
\begin{array}{cl}
\min\limits_{x\in X}  & x_3  \\
\st &  \begin{bmatrix}x_3-5x^2_1-5x^2_2+u^2_1+u^2_2-x_1(-u_1+u_2+5)-x_2(u_1-u_2+3)\\
	0.2-x^2_1-u^2_1\\
	0.1-x^2_2-u^2_2
\end{bmatrix}
\geq 0 \\
& \forall u\in U=[-0.2,0.2]^2,\\
&x\in X=[-100,100]^3.\\
\end{array}
\right.
\end{equation*}
This SIP is solved by Algorithm~\ref{algorithm} with  $3$ loops.
The computational results are shown in
Table~\ref{Computational results for SIP_ex_1}.
\begin{center}
\begin{longtable}{cll}
\caption{Computational results for Example~\ref{SIP_ex_1}}
\label{Computational results for SIP_ex_1}\\
\specialrule{.1em}{0em}{0.1em}  	
$k$ & \hspace{4em} $(\hat{x}_k,\hat{u}_k)$ & function values \\ %[0.5ex]
\midrule
$0$ & $\hat{x}_0=(  0.0000,   0.0000, -100.0000)$ & $ f_0 =  -100.0000$\\
&$\hat{u}_0=(0.0000 ,   0.0000 )$ &
$ \hat{g}_0=  -100.0000 $ \\ \cmidrule{1-3}
$1$ & $\hat{x}_1=(  -0.4000,   -0.2449,   -1.6348)$ & $ f_1 =  -1.6348$\\
&$\hat{u}_1=( 0.0775 ,  -0.0775 )$ &
$ \hat{g}_1=  -0.0120 $ \\
\cmidrule{1-3}
$2$ & $\hat{x}_2=(-0.4000,   -0.2449,   -1.6228)$ & $ f_2 =  -1.6228$\\
&$\hat{u}_2=( -0.2000 ,   0.2000)$ &
$ \hat{g}_2=  -7.0\cdot 10^{-9} $ \\ \cmidrule{1-3}
&Output: $x^*=\hat{x}_2$, $f^*=f_2$ & \text{time:} \, $   1.62s$  \\
\specialrule{.1em}{0em}{0.1em}
\end{longtable}
\end{center}
\end{eg}

\begin{eg}\label{SIP_ex_2}
Consider the following SIP
\begin{equation*}
\left\{
\begin{array}{cl}
\min&\left(x_{1}-x_{2}\right)\left(x_{3}-x_{4}\right)+
\left(x_{1}-x_{3}\right)
\left(x_{2}-x_{4}\right)+
x_{1} x_{2}-x_{2} x_{3}+x_{3} x_{4}  	\\
\st & \begin{bmatrix}
	-2x_1x_2u_2u_3+u^2_1x_4+x_1u_2-x_2u_3+x_3u_1+x_4+1\\
	2u_1x_2-u_3x_4+2u_2u_3+2x_1x_4-x_2x_3+1\\
\end{bmatrix}\geq 0\quad\forall u\in U,\\
&x\in X,
\end{array}
\right.
\end{equation*}
where the constraining sets are
\[
\begin{array}{rcl}
X &=& \left\{x\in\mathbb{R}^4\,\middle\vert\,\begin{array}{ll}
	10-x^Tx\geq 0,\,3x_1x_2-2x^2_3+4x_4 \geq 0\\
	x_1x_2x_3-1\geq 0,\,-3x_3x_4+x_2+2x_3\geq 0\\
\end{array}\right\},
    \\
U &=& \left\{u\in\mathbb{R}^3\,\middle\vert\,\begin{array}{ll}
	u_1u_2u_3-1\geq 0,\\
	-2u^2_3+u_1+u_2+3\geq 0,\\
	-u^Tu+u_1+u_3+5\geq 0,
\end{array}\right\}.
\end{array}
\]
This SIP is solved by Algorithm~\ref{algorithm} with  $2$ loops.
The computational results are shown in
Table~\ref{Computational results for SIP_ex_2}.
\begin{center}
\begin{longtable}{c l l }
\caption{Computational results for Example~\ref{SIP_ex_2}}
\label{Computational results for SIP_ex_2}\\
\specialrule{.1em}{0em}{0.1em}  	
$k$ & \hspace{4em} $(\hat{x}_k,\hat{u}_k)$ & function values \\
\midrule
$0$ & $\hat{x}_0=(1.2517,   -1.3709  , -1.3383   , 2.1824  )$ & $f_0= -24.9074$\\
&$\hat{u}_0=(2.5046 ,   0.2357 ,   1.6941)
$ & $\hat{g}_0= -5.1372	$ \\  \midrule
$1$ & $\hat{x}_1=(1.7887 ,  -0.9005 ,  -1.3106 ,   2.0669  )$ & $f_1= -23.7793$\\
&$\hat{u}_1=(2.5046  ,  0.2357 ,   1.6941)
$ & $\hat{g}_1=   3.4\cdot 10^{-8}$ \\ \cmidrule{1-3}
&Output: $x^*=\hat{x}_2$, $f^*=f_2$ & \text{time:} \, $   1.22 s$  \\
\specialrule{.1em}{0em}{0.1em}
\end{longtable}
\end{center}
\end{eg}

\subsection{GSIP examples}

This subsection reports numerical examples for GSIPs.

\begin{eg}\label{Box change variable}
Consider the GSIP from Example~\ref{Alexander 2}.
Make the substitution $y\coloneqq u^5$, then Example~\ref{Alexander 2} becomes:
\[
\left\{
\begin{array}{cl}
	\min & -x_1\\
	\st  & y-3x^2_2\geq 0\quad\forall y\in Y(x),\\
	& x\in X,
\end{array}
\right.
\]
where the constraining sets are
\[
X=[0,1]^2, \quad
Y(x)=\left\{y\in\mathbb{R} :
\begin{array}{cl}
	-32\leq y\leq 0, 1-4x^2_1-x^2_2\leq y
\end{array}\right\}.
\]
For each $x\in X$, note that $1-4x^2_1-x^2_2 \ge  -32$.
When $1-4x^2_1-x^2_2 > 0$, the set $Y(x)$ is empty.
We solve this GSIP in two cases:

\noindent
Case I: If $x\in X\cap \{1-4x^2_1-x^2_2\leq 0\}$, then
$Y(x)=\{1-4x^2_1-x^2_2\leq y\leq 0 \}$, which is a box constraint.
The polynomial extension can be chosen as in \reff{PE:box}.
For this case of the GSIP, Algorithm~\ref{algorithm} terminates in $2$ loops
and returns $x^*\approx(  0.5000, \,   0.0000), f^*\approx  -0.5000$.
The computational results are shown in
Table~\ref{result for Box change variable}.

\noindent
Case II: If $x\in X\cap \{1-4x^2_1-x^2_2> 0\}$, then $Y(x)=\emptyset $
and the GSIP is equivalent to
\be  \label{Example 6.3 Case 2}
\left\{
\begin{array}{cl}
	\min \, &  -x_1\\
	\st  \, & x\in X\cap \{1-4x^2_1-x^2_2> 0\}.
\end{array}
\right.
\ee
When this polynomial optimization is solved, the strict inequality
$1-4x^2_1-x^2_2> 0$ is treated as $1-4x^2_1-x^2_2\geq  0$.
Solving \reff{Example 6.3 Case 2} gives
$\tilde{x} \approx  (0.5000 ,  0.0001), f(\tilde{x}) \approx   -0.5000$.
The minimizer $\tilde{x}$ is the same as the one for Case I.
The constraint $1-4x^2_1-x^2_2 \geq  0$ is active at $\tilde{x}$.
By comparison with Case I, $\tilde{x}$ is a global minimizer for the GSIP.
\begin{center}
\begin{longtable}{c l l }
\caption{Computational results for Example~\ref{Box change variable}~Case~I.}
\label{result for Box change variable}\\
\specialrule{.1em}{0em}{0.1em}  	
$k$ & \hspace{2em} $(\hat{x}_k,\hat{y}_k)$ & function values \\ \midrule
$0$ & $\hat{x}_0 = (1.0000,    0.3251)$ & $f_0=   -1.0000$ \\
& $\hat{y}_0 = -3.1057$ & $\hat{g}_0 =  -3.4229$ \\ \midrule
$1$ & $\hat{x}_1 = ( 0.5000,    0.0000)$ & $f_1=  -0.5000$ \\
& $\hat{y}_1 = 0.0000$ & $\hat{g}_1 =   -1.3\cdot 10^{-9}$ \\ \cmidrule{1-3}
& \text{Output:} $x^*=\hat{x}_1$, $f^*=f_1$, $g^*=g_1$ & \text{time:} $0.68s$ \\
\specialrule{.1em}{0em}{0.1em}
\end{longtable}
\end{center}
\end{eg}

\begin{eg}\label{GSIP_ex_9}
Consider the following GSIP
\begin{equation*}
\left\{
\begin{array}{cl}
	\min & x^2_1x_2-x^2_2x_3-x_1x_2x_3+(x_1+1)^2\\
	\st  & \begin{bmatrix}
		x^2_1u_1+x_2u^2_2-x_3u^2_1-x_1x_3u_2+1\\
		x_1x_2u_1u_2+x_1x_3(u_1+u_2)-x_3u_2-1
	\end{bmatrix}\geq 0\quad\forall u \in U(x),\\
	& x\in X,
\end{array}
\right.
\end{equation*}
where the constraining sets:
\[
\begin{array}{rcl}
X &=& \left\{x\in\mathbb{R}^3\,\middle\vert\,\begin{array}{ll}
	25-x^Tx\geq 0,x_1-x^2_2\geq 0,\\
	x_3-x_2-1\geq 0\\
\end{array}\right\},
\\
U(x) &=& \left\{u\in\mathbb{R}^2\,\middle\vert\,\begin{array}{ll}
	0\leq u_1\leq x_1,\\
	x_2\leq u_2\leq x_3
\end{array}\right\}.
\end{array}
\]
For all $x\in X$, we have $0\leq x_1,\,x_2<x_3$, so $U(x)$
is nonempty and has box constraints. The polynomial extension
is given in Section~\ref{sec:Poly}.	However,
in the $2$nd loop  at Algorithm~\ref{algorithm},
the subproblem $(P_2)$ is infeasible, so Example~\ref{GSIP_ex_9} is infeasible.
The result is shown in Table~\ref{Computational results for GSIP_ex_9}.
\begin{center}
\begin{longtable}{c l l }
\caption{Computational results for Example~\ref{GSIP_ex_9}}
\label{Computational results for GSIP_ex_9}\\
\specialrule{.1em}{0em}{0.1em}
$k$ & \hspace{5em} $(\hat{x}_k,\hat{u}_k)$ & function values \\ \midrule
$0$ & $\hat{x}_0=(  4.4034,  -2.0984,   -1.0984)$ &  $f_0 = -16.8047$  \\
& $\hat{u}_0=(0.0000,   -2.0984)$  &  $\hat{g}_0 =   -18.3900$ \\ \midrule
$1$ & $\hat{x}_1=(4.4449,   -2.1083,    0.8933)$ &  $f_1 = -7.6063$  \\
& $\hat{u}_1=( 4.4449,    0.8933)
$  &  $\hat{g}_1 = -17.8112$ \\ \cmidrule{1-3}
$2$ &\textit{Moment relaxation is infeasible} \\
\cmidrule{1-3}
& \textrm{Output}: The GSIP in Example~\ref{GSIP_ex_9} is infeasible.  \\
\specialrule{.1em}{0em}{0.1em}
\end{longtable}
\end{center}
\end{eg}

\begin{eg}\label{GSIP ex ellipse case}
Consider the following GSIP
\[
\left\{
\begin{array}{cl}
\min & x_{1}^{2} x_{2}-x_{1} x_{2}^{2}-
x_{1} x_{2} x_{3}+x_{3} x_{4} x_{5}+x_{3}^{3}\\
\st  & \begin{bmatrix}
	x_1u_1+x_2x_5u_5-x_3u_2u_3-x_4u_4-1\\
	x_1x_2u_1u_2-x_3x_4u_1-x_5u_3u_4u_5-5\\
	u_1u_2+x_5u_3-(x_1-u_1)^2-(x_3-u_4)^2-u^2_5+1
\end{bmatrix}\geq 0\quad\forall u\in U(x)\\
&x\in X,
\end{array}
\right.
\]
where the constraining sets are
\[
\begin{array}{rcl}
X &=& \left\{x\in\mathbb{R}^5\,\middle\vert
\begin{array}{ll}
	25-x^Tx\geq 0,
	x_1x_2-1\geq 0, x_3\geq 0,\\ -x_2+x_3x_5-3\geq 0,x_3x_4-2\geq 0 \\
\end{array}
\right\},
 \\
U(x) &=& \left\{u\in\mathbb{R}^5:
\begin{array}{ll}
	\|u-(x_1,\,x_2,\,0,\,0,\,0)\| \le \sqrt{3}x_3
\end{array}
\right\}.
\end{array}
\]
Since $x_3\geq 0$ for all $x\in X$, the set $U(x)$ has ball constraints,
so we can choose the polynomial extension
\[
q(x) =  \frac{1}{ \hat{x}_3 }
\begin{bmatrix}
	x_1\hat{x}_3+(u_1-x_1)x_3,\,
	x_2\hat{x}_3+(u_2-x_2)x_3,\,
	\hat{u}_3x_3,\,
	\hat{u}_4x_3,\,
	\hat{u}_5x_3
\end{bmatrix} .
\]
This GSIP is solved by Algorithm~\ref{algorithm} with $4$ loops.
The computational results are shown in
Table~\ref{result for GSIP ex ellipse case}.
\begin{center}
\begin{longtable}{c l l }
\caption{Computational results for Example~\ref{GSIP ex ellipse case}}
\label{result for GSIP ex ellipse case}\\
\specialrule{.1em}{0em}{0.1em}  	
$k$ & \hspace{4em} $(\hat{x}_k,\hat{u}_k)$ & function values \\ \midrule
$0$ & $\hat{x}_0=(  -4.0182,   -2.1036,    1.5910,    1.2571,    0.5634)$ & $f_0=    -24.4776$\\
&$\hat{u}_0=(  -4.2767,    0.6399,   0.0000,   0.0000 ,   0.0000)$ &
$\hat{g}_0= -19.5772$ \\
\midrule
$1$ & $\hat{x}_1=(-3.7475,   -2.3203,    0.9845,    2.0314,    0.6903)$ &
$f_1= -18.6352  $ \\
&$\hat{u}_1=(  -2.8124,  -1.5519,   -0.1886,   -1.1863,  0.0000)$ &
$\hat{g}_1 =  -0.3528$ \\ \midrule
$2$ & $\hat{x}_2=( -3.7162,   -2.3348,    0.9605,    2.0823,    0.6926) $ &
$f_2=  -18.0488 $\\
&$\hat{u}_2=(-2.7811,   -1.5813,   -0.1876,   -1.1358,    0.0000)$ &
$\hat{g}_2= -0.0010$ \\ \midrule
$3$ & $\hat{x}_3=(-3.7163,   -2.3344,    0.9603,    2.0827 ,   0.6931	)$ & $f_3= -18.0471$\\
&$\hat{u}_3=( -2.7814,   -1.5808,   -0.1878,   -1.1356,   0.0000)$ &
$\hat{g}_3=    -1.8\cdot 10^{-7}$ \\ \cmidrule{1-3}
& \text{Output:} $x^*=\hat{x}_3$, $f^*=f_3$	& \text{time:} $2.01s$ \\
\specialrule{.1em}{0em}{0.1em}   	
\end{longtable}
\end{center}
\end{eg}

\begin{eg}\label{simplex constraint case}
Consider the following GSIP $(\text{let} \,\, x_0:=1)$:
\[
\left\{
\begin{array}{rl}
\min&\sum\limits_{i=0}^{3}\prod\limits_{j\not=i}(x_i-x_j) \\
\st  &\begin{bmatrix}
	(x_1x_2-u_1)(x_3-u_2)+x_2u_1u_2+x_1x_3+2\\
	2x_1u_1-x_3u_1u_2+x_2+1
\end{bmatrix}\geq 0\quad\forall
u\in U(x),\\
&x\in X,
\end{array}\right.
\]
where the constraining sets
\[
X=\left\{x\in\mathbb{R}^3\,\middle\vert\,
\begin{aligned}
	&x_1+x_2-x_3\geq 0, \, x_1-x_2\geq 0,\, 1-x_1-x_3\geq 0,\\
	&2-2x_1+x_2\geq 0,\, 3x_2+x_3-1\geq 0
\end{aligned}
\right\},
\]
\[
U(x)=\left\{u\in\mathbb{R}^2\,\middle\vert\begin{array}{ll}
	e^Tx\geq e^Tu,\,x_1u_2-u_1-x_3-x_1x_2\geq 0,\,\\
	x_1-u_2+x_2u_1\geq 0,\,-u_1\geq 0
\end{array}
\right\}.
\]
The polynomial extension $q $ is obtained
by solving the system \reff{polyext:q:num}.
This GSIP is solved by Algorithm~\ref{algorithm} with $4$ loops.
The computational results are shown in
Table~\ref{result for simplex constraint case}.
\begin{center}
\begin{longtable}{c l l}
\caption{Computational results for Example~\ref{simplex constraint case}}
\label{result for simplex constraint case}\\
\specialrule{.1em}{0em}{0.1em}  	
$k$ & \hspace{5em} $(\hat{x}_k,\hat{u}_k)$ & function values \\ %[0.5ex]
\midrule
$0$ & $\hat{x}_0=(2.0000,    2.0000,   -5.0000)$ & $f_0=-288.0000$\\
& $\hat{u}_0=(  -1.0000, 0.0000)$&$\hat{g}_0=-33.0000$ \\
\midrule
$1$ & $\hat{x}_1=(0.7195  ,  0.7195,   -1.1585)$ & $f_1=-7.4436$\\
& $\hat{u}_1=(   -0.8177 ,   0.1312)$&$\hat{g}_1=-0.6331$ \\
\midrule
$2$ & $\hat{x}_2=(0.6496 ,   0.6496  , -0.9489)$ & $f_2=-4.7405$\\
& $\hat{u}_2=(  -0.8148 ,   0.1203)$&$\hat{g}_2=-0.0025$ \\
\midrule
$3$ & $\hat{x}_3=(0.6493  ,  0.6493,   -0.9480)$ & $f_3= -4.7306$\\
& $\hat{u}_3=(  -0.8148,    0.1203)$&$\hat{g}_3= 2.5\cdot 10^{-8}$ \\
\cmidrule{1-3}
& \text{Output:} $x^*=\hat{x}_3$, $f^*=f_3$, $g^*= \hat{g}_3$ & \text{time:} $
6.53s$ \\
\specialrule{.1em}{0em}{0.1em}
\end{longtable}
\end{center}
\end{eg}

\begin{eg}\label{four_application_examples}
We show how to find maximum ellipsoid contained inside
a region. This can be done by solving the GSIP
\eqref{design centering problem}. The ellipsoid can be represented as
\[
U(x) \, = \,  \left\{u\in\mathbb{R}^2:
\begin{array}{ll}
	(u - a)^T P(x)^{-1} (u - a) \leq 1
\end{array}
\right\}.
\]
The center of the ellipsoid is $a = (x_1, x_2)$ and the shape matrix is
\[
P(x) = \bbm x_3 & 0 \\ x_4 & x_5 \ebm
\bbm x_3 & x_4 \\ 0 & x_5 \ebm .
\]
The GSIP for finding the maximum ellipsoid is
\begin{equation*}
\left\{\begin{array}{cl}
	\min & -\pi\det\begin{bmatrix}x_3&x_4\\0&x_5\end{bmatrix} \\
	\st  & 	g(u)\geq 0\quad \forall u\in U(x),\\
	& x\in X,
\end{array} \right.
\end{equation*}
where $X= \{ x \in \re^5 : \, 100-x^Tx\geq 0, x_3\geq 0, x_5\geq 0 \}$ and
\[
g(u) = (4u_1+u^2_2-4,\, -u_1+2u_2+4,\, 8-(u_1-1)^2-u^2_2,\, 3-u_1u_2 ) \ge 0 .
\]
We select the polynomial extension as in \reff{polyext:Ellipsoid}.
After $9$ loops, we obtained the optimizer
\[
x^* \approx (1.5084, 1.0587, 1.4203, -1.4097, 0.9039 ) .
\]
The result is shown in Figure~\ref{SetFigures}.
It took around $14$ seconds.
\begin{figure}[htb]
\centering
\begin{tabular}{c}
	\includegraphics[width=.45\textwidth]{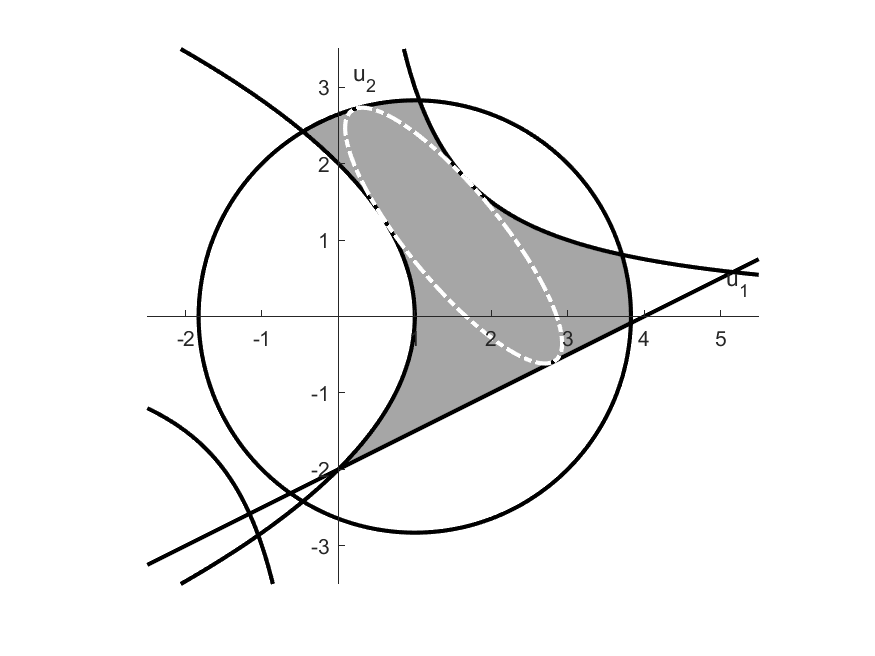}
\end{tabular}
\caption{The shaded area is the region of $g(u) \ge 0$
	in Example~\ref{four_application_examples}.
	The dotted line is the boundary of
	the maximum ellipsoid inscribed in the shade. }
\label{SetFigures}
\end{figure}
\end{eg}

\subsection{GSIPs with convex infinite constraints}
\label{ssc:cvxgsip}

\begin{eg}\label{LME_example1}
Consider the following GSIP from \cite{jungen2023libdips}:
\begin{equation*}
\left\{
\begin{array}{cl}
	\min& x^2_1+x^2_2\\
	\st &  g(x,u)= (u_1-x_1)^2+(u_2-x_2)^2-1\geq 0\quad  \forall u\in U(x),\\
	&x\in \mathbb{R}^2,\\
\end{array}
\right.
\end{equation*}
where
$
U(x) = \{u\in\mathbb{R}^2 :
u_1-x_1\geq 0,\, u_2\geq 0\}.	
$
Both the function $g(x,u)$ and the set $U(x)$ are convex in $u$.
Since $U(x)$ is given by linear inequality constraints,
the KKT condition \reff{cvxKKT:mingjxu} holds.
Solve the polynomial optimization problem~\reff{cvx:GSIP} with variables $x,z,\lambda$,
where $z$ served as the parameter $u$ in $g(x,u)$,
and $\lambda$ is the Lagrange multiplier.
The computed optimal solutions and values are
\[
\begin{aligned}
	&x^*\approx(  6.2\cdot 10^{-6}, -1.0000), \,
	z^*\approx (1.3\cdot 10^{-4}, 5.1\cdot 10^{-7}), \\
	& \lambda^*\approx(2.3\cdot 10^{-4},   2.0000),\,
	f(x^*)\approx   1.0000,\,
	g(x^*, z^*) \approx  -1.6\cdot 10^{-5} .
\end{aligned}
\]
It took around $3.7$ seconds.
\end{eg}

\begin{eg}  \label{LME_example2}
Consider the following GSIP example:
\begin{equation*}
	\left\{
	\begin{array}{cl}
		\min& -x_1x_2+x_1+2x_2\\
		\st &  g(x,u)= \begin{bmatrix}
			\displaystyle u_1+u_2+\frac{x^2_2}{x_1+u_1}+\frac{x^2_1}{x_2+u_2}-2\\
			x_1u_2+\displaystyle\frac{x^2_2+u^2_2}{1+u_1}-x^2_2
		\end{bmatrix}\geq 0\quad  \forall u\in U(x),\\
		&x\in X,\\
	\end{array}
	\right.
\end{equation*}
where the constraining sets are
\[
\baray{rcl}
	X &=& \big \{x\in\mathbb{R}^2: 4-x^Tx\geq 0,\,x_1x_2\geq
	\frac{1}{2},\,x_1\geq 0 \big \}, \\
	U(x) &=& \{u\in\mathbb{R}^2:
	u_1\geq 0,\, u_2-x_2u_1\geq 0,\, 2x_1-x_2u_1-u_2\geq 0\}.
\earay
\]
The constraining function $g_1(x,u),\,g_2(x,u)$ are both convex in $u=(u_1,\,u_2)$
over $\mathbb{R}_{+}\times \mathbb{R}$, and the set $U(x)$ is convex in $u$.
For any $x\in X$, we have $x>0$, and $\nabla_u h\not=0$.
Consider point $\hat{u}=(\frac{x_1}{2x_2},\,x_1)$, then
$h(x,\hat{u})=(\frac{x_1}{2x_2},\,\frac{x_1}{2},\,\frac{x_1}{2})> 0$,
so the KKT condition holds for $u\in U(x)$.
The matrix polynomial $T(x,u)$ satisfying equation \reff{RLME equation} is:
\[T(x,u)=
\begin{bmatrix}
	2x_1-2x_2u_1& 2x_1x_2-2x_2u_2 & 2x_2& 2x_2& 2x_2\\
	-u_1 & 2x_1-u_2 & 1& 1& 1\\
	-u_1 & -u_2 & 1& 1& 1\\
\end{bmatrix},
\]
for the denominator vector $\phi(x,u)=2x_1.$
Solving the optimization~\reff{cvx:GSIP:lme} with variables $x,z_1,z_2$,
where $z_i$ served as the parameter $u$ in $g_i(x,u),\,i=1,2$,
we get the optimizers and optimal values:
\[
\begin{aligned}
	&x^*\approx(  1.1348,    0.4406), \,z^*_1\approx ( 9.8\cdot 10^{-9},    0.6941),
	\,z^*_2\approx(1.9,    2.1)\cdot 10^{-9},\\
	&f(x^*) \approx   1.5160,\, \big( g_1(x^*, z^*_1), g_2(x^*, z^*_2) \big)
	\approx(4.9, 0.2)\cdot 10^{-7}. %%  \text{time}: 129.52 s.
\end{aligned}
\]
It took around $129$ seconds.
This GSIP can also be solved by Algorithm~\ref{algorithm}.
It terminates within $3$ iterations and returns the same results.
Algorithm~\ref{algorithm} took around $6$ seconds.
This is faster than solving  \reff{cvx:GSIP:lme},
because \reff{cvx:GSIP:lme} has extra variables $z_1,z_2$,
which made the resulting Moment-SOS relaxations more expensive to solve.
\end{eg}

\subsection{GSIPs from existing references}		
\label{ssc:appendix}

This subsection gives two tables of computational results
for SIP and GSIP examples from earlier existing references.
All of these examples are successfully solved by  Algorithm~\ref{algorithm}.
The SIP examples are listed in Appendix~\ref{AppendixA}
and the computational results are reported in Table~\ref{all result for appendix A}.
The GSIP examples are listed in Appendix~\ref{AppendixB}
and their computational results are reported in Table~\ref{all result for appendix B}.
In these examples, if the constraining set $X$
is not given in the original references,
we set it to be $X=[-100,100]^n$.

\begin{center}
\begin{longtable}{ccll}
\caption{Computational results for SIPs in Appendix~\ref{AppendixA}}
\label{all result for appendix A}\\
\specialrule{.1em}{0em}{0.1em}
Problem & \# of  & $\hspace{2em} (x^*,u^*)$ & $(f^*,g^*,\text{time})$ \\
& loops &  &  \\
\midrule
$\eqref{Watson_1}$ & $2$ & $x^* = (-0.7500, -0.6180)$ & $f^* = 0.1945$ \\
& &	$u^*= 4.1\cdot 10^{-6}$ & $g^* = -3.9\cdot 10^{-9}$ \\
& & & \text{time:}  $0.70s$ \\
\cmidrule{1-4}
$\eqref{Waston 2}$& $3$& $x^* = (-1.0000,0.0000,0.0000)$ & $f^* = 1.0000$ \\
& &	$u^* = (0.0, 0.3)\cdot 10^{-5}$  & $g^* = -8.2\cdot 10^{-10}$ \\
& &  & \text{time:} $0.98s$ \\
\cmidrule{1-4}
$\eqref{Wang_3}$& $2$& $x^*=(-0.3,1.7)\cdot 10^{-11}$ & $f^* = 1.7\cdot 10^{-11}$ \\
& &	$u^* =  5.7\cdot 10^{-11}$ &  $g^* = 3.1\cdot 10^{-11}$  \\
& &  &  \text{time:}  $0.69s$ \\
\cmidrule{1-4}
$\eqref{Wang_PMI} $& $5$& $x^*= (-1.2808,-1.2808)$  & $ f^* = -2.5616 $  \\
&  & $u^* = (-0.9294,   -0.2610,   -0.2610)$ &  $g^* = 7.7\cdot 10^{-10}$ \\
&  &  &  \text{time:} $1.29s$  \\	
\cmidrule{1-4}
$\eqref{Lemonidis_3}$& $11$& $ x^*=(  3.0000,0.0000,0.0000$ &  $f^*= -12.0000$ \\
& & \hspace{3em} $		-0.0019,	0.0000, 0.0019 )$ &  $g^* = -9.4\cdot 10^{-7}$  \\
& & $u^* = ( -0.0298,    0.0086)$ &  \text{time:} $3.40s $\\
\cmidrule{1-4}
$\eqref{Teo 1} $& $2$& $x^* = (0.7071, 0.7071)$  &  $f^* = 0.3431 $  \\
& &	$u^* = ( 0.7853, 0.7853)$ & $g^* =  2.4\cdot 10^{-9}$ \\
& &  & \text{time:} $0.91s$\\
\cmidrule{1-4}
$\eqref{betro_sum} $ & $13$
& $x^*=( 0.5000,    0.2501,    0.1226,    0.0798,$ & $f^* = 0.6931$  \\
& & \hspace{2.5em}$ -0.0299,    0.1297,   -0.1026,    0.0502
)$ & $g^* = -4.1\cdot 10^{-7}$ \\
& & $u^* = 0.1637 $ &  \text{time:} $2.06s$ \\
\cmidrule{1-4}
$\eqref{Flouda 1} $ & $ 7$ & $x^*=(-0.0280, 4.0001, -4.0002, 0.0280)$ & $f^* =  0.0280  $  \\
& &$u^*=0.1502$ &  $g^* =
-1.3 \cdot 10^{-7}$ \\
& & & \text{time:} $2.96s$  \\
\cmidrule{1-4}
$\eqref{Zakovic-minmax} $& $ 9 $ & $x^*= ( 1.6953, 0.0000
)$ & $f^* =  1.4039$  \\
&	& $u^* = ( 0.7185, 0.0000)$ & $ g^* =   -1.4\cdot 10^{-7} $ \\
&  &  &	\text{time:} $2.52s$ \\
\cmidrule{1-4}
$\eqref{lemondis_paper}$ & $2$ & $x^* = ( -0.0001,    0.4999)$ &  $f^* =  -0.2500$  \\
&  &  $u^*= 3.5\cdot 10^{-12} $ & $ g^* =  -1.0\cdot 10^{-8} $  \\
&  &  & \text{time:}  $0.59s$\\
\specialrule{.1em}{0em}{0.1em}
\end{longtable}
\end{center}

\begin{center}
\begin{longtable}{ccll}
\caption{Computational results for Appendix~\ref{AppendixB}}
\label{all result for appendix B}\\
\specialrule{.1em}{0em}{0.1em}
Problem & \# of  & $\hspace{2em} (x^*,u^*)$ & $(f^*,g^*,\text{time})$ \\
& loops &  &  \\
\midrule
$\eqref{Aboussoror1}$& $2$ & $x^*= 0.5000$ &  $f^* = 1.7500$ \\
&  & $u^* = (-0.7500 ,  -0.5000)$ &  $g^* = 3.1\cdot 10^{-10}$  \\
&  &   &  \text{time:}  $0.58s$\\
\cmidrule{1-4}
$\eqref{Alexander 1}$ & $2$ & $ x^* = (0.9999 ,   0.9999)$ & $f^* = -0.5000$ \\
&  & $u^* =  ( 0.7071,    0.7071,    0.0000
)$  & $g^* = 2.8\cdot 10^{-8}$ \\
&  &   & \text{time:}  $  0.82s$  \\
\cmidrule{1-4}
$\eqref{M.Diehl 1}$ & $1$ & $x^*=( 0.0000, 1.0000)$  &  $f^* = 5.8284$   \\
& & $u^*= 9.0\cdot 10^{-10}$ & $g^* =  9.0\cdot10^{-10}$  \\
& & & \text{time:} $0.53s$   \\
\cmidrule{1-4}
$\eqref{jongen 1}$ & $2$ & $x^* = (0.0000,0.0000)$ &
$f^* =-1.6\cdot 10^{-11}$  \\
& & $u^* = -1.1\cdot 10^{-4}$ & $g^* =    -2.0\cdot 10^{-11}$  \\
& &  & \text{time:}  $0.43 s$          \\
\cmidrule{1-4}
$\eqref{Alexander_4}$& $3$& $x^*=(-0.6180,0.0000,0.0000)$ &  $f^* = 0.3820$  \\
&	& $u^*=  (0.6180,0.0000)$ &  $g^* =  3.1\cdot 10^{-9}$  \\
&   &  &  \text{time:} $1.09s$  \\
\cmidrule{1-4}
$\eqref{Alexander_5}$ & $2$ & $x^*=( -0.5000,-0.5000, 0.0000)$ &  $f^* = 0.5000$  \\
&  &  $u^* =   4.1\cdot 10^{-11}$ &  $g^* = -9.7\cdot 10^{-7} $ \\
&  &   & \text{time:}  $ 1.54s$  \\
\cmidrule{1-4}
$\eqref{Alexander_6}$ & $2$ & $x^*=  ( -1.0000, 0.0000)$ &  $f^* = -1.0000$   \\
& & $x^*=  (0.0000, -1.0000)$ &   $g^* =  7.5\cdot 10^{-9} $   \\
&  &  $u* = 5.8 \cdot 10^{-10}$   &   \text{time:} $1.94s$  \\
\cmidrule{1-4}
$\eqref{Alexander_7}$ & & $x^* = (-1.0000,0.2500,0.2500)$ & $f^* = 2.9361$  \\
&  &    &   \text{time:}  $ 0.55 s$  \\
\cmidrule{1-4}
$\eqref{Harwood 1}$ & $4$  & $x^* = (-1.0000,-1.0000,1.0000,-0.1510)$ & $f^* =  -1.6980$  \\
&  & $u^* = (0.7380, -0.1510)$ &  $g^* =  -2.3\cdot 10^{-7}$  \\
&  &    & \text{time:} $203.69s$\\
\cmidrule{1-4}
$\eqref{O.Stein 1}$ & $ 10$ & $x^* = (  2.0125,   -0.4997,    2.2164,    0.5003	)$ &
    $f^* = -3.4838$ \\
	& &  $u^*=( 3.6575,   -0.1644)$ & $g^* = -1.1 \cdot 10^{-7}$ \\
	&  & 	&  \text{time:} $2.86s$\\
	\specialrule{.1em}{0em}{0.1em}
\end{longtable}
\end{center}

\section{Conclusions and Discussions}
\label{sec:COnclu}

This paper proposes a hierarchy of polynomial optimization relaxations
to solve GSIPs. The optimization problems are not assumed to be convex.
When the constraining set $U(x)$ depends on $x$,
we use polynomial extensions to construct polynomial optimization relaxations.
We have also demonstrated how to get polynomial extensions for various feasible sets.
Moment-SOS relaxations are applied to solve the polynomial optimization.
The method is given in Algorithm~\ref{algorithm}.
Its efficiency is demonstrated by extensive examples.
In particular, the SIP can be solved as a special case of the GSIP.
An interesting fact is that the finite convergence
of Algorithm~\ref{algorithm} is observed numerically
in all our computational examples.

There is much interesting future work to do.
For instance, how do we get polynomial extensions for general set $U(x)$?
How do we efficiently solve the polynomial relaxation $(P_k)$
when there are a lot of constraints?
What are conveniently checkable conditions for
Algorithm~\ref{algorithm} to have finite convergence?
These questions are mostly open, to the best of the authors's knowledge.

\bigskip
\bigskip \noindent
{\bf Acknowledgements}   
The research is partially supported by the NSF grant DMS-2110780.

\appendix
\section{The SIP examples from references}
\label{AppendixA}

The SIPs in this appendix are from existing work.
Their computational results are shown in Table~\ref{all result for appendix A}.

\begin{eg}\label{Watson_1}
The SIP from \cite{coope1985projected,wang2015feasible}:
\begin{equation*}
\left\{
\begin{array}{cl}
	\min& \frac{x^2_1}{3}+\frac{x_1}{2}+x^2_2\\
	\st & g(x,u)=-(1-x^2_1u^2)^2+x_1u^2+x^2_2-x_2\geq 0\quad\forall u\in U,\\
	&U=[0,1].
\end{array}
\right.
\end{equation*}
\end{eg}

\begin{eg}\label{Waston 2}
The SIP from \cite{coope1985projected}
\begin{equation*}
\left\{
\begin{array}{cl}
	\min& x^2_1+x^2_2+x^2_3 \\
	\st & g(x,u)=-x_1(u_1+u^2_2+1)-x_2(u_1u_2-u^2_2)\\
	&\quad\quad\quad\quad-x_3(u_1u_2+u^2_2+u_2)-1\geq 0\quad\forall u\in U,\\
	&U=[0,1]^2.
\end{array}
\right.
\end{equation*}
\end{eg}

\begin{eg}\label{Wang_3}
The SIP from \cite{wangliguofeng2014semidefinite}
\begin{equation*}
\left\{
\begin{array}{cl}
	\min& x_2 \\
	\st & g(x,u)=-2x^2_1u^2+u^4-x^2_1+x_2\geq 0\quad\forall u\in U,\\
	&U=[-1,1].
\end{array}
\right.
\end{equation*}
\end{eg}

\begin{eg}\label{Wang_PMI}
The SIP from \cite{wangliguofeng2014semidefinite}
\begin{equation*}
\left\{
\begin{array}{cl}
	\min& x_1+x_2 \\
	\st & G(x)= \begin{bmatrix}
		4-x_{1}^{2}-x_{2}^{2} & x_{1} & x_{2}\\
		x_{1} & x_{2}^{2}-x_{1} & x_{1} x_{2} \\
		x_{2} & x_{1} x_{2} & x_{1}^{2}-x_{2}
	\end{bmatrix}\succeq 0
\end{array}
\right.
\end{equation*}
For the PSD matrix constraint, we can convert $G(x)\succeq 0$ to the equivalent form
$g(x,u)=u^TG(x)u \geq 0$ with $U  = \{ u\in\mathbb{R}^3 :  u^Tu=1\}$.
\end{eg}

\begin{eg}\label{Lemonidis_3}	
The SIP from \cite{lemonidis2008global}
\begin{equation*}
\left\{
\begin{array}{cl}
	\min & -4x_1-\frac{2}{3}(x_4+x_6) \\
	\st & g(x,u)= -x_1-x_2u_1-x_3u_2-x_4u^2_1-x_5u_1u_2-x_6u^2_2\\
	&\quad\quad\quad\quad+(u_1-u_2)^2(u_1+u_2)^2 +3 \geq
	0\quad\forall u\in U,\\
	&U=[-1,1]^2.
\end{array}
\right.
\end{equation*}
\end{eg}

\begin{eg}\label{Teo 1} 	
The SIP from \cite{teo2000computational}
\begin{equation*}
\left\{
\begin{array}{cl}
	\min& (x_1+x_2-2)^2+(x_1-x_2)^2+30 \big( \min \{0, x_1-x_2 \} \big)^2 \\
	\st & g(x,u)=-x_1\cos(u)-x_2\sin(u)+1\geq 0\quad\forall u\in U, \\
	&U=[0,\pi].
\end{array}
\right.
\end{equation*}
We introduce the new variable $x_3  \coloneqq  \min \{0, x_1-x_2 \}$,
which is the same as \[ 2x_3 =x_1-x_2 - |x_1-x_2| \]
and is equivalent to
\[
x_1 - x_2 - 2x_3 \ge 0, \quad (x_1 - x_2 - 2x_3)^2 = (x_1 - x_2 )^2.
\]
Let $w_1=\cos(u)$ and $w_2=\sin(u)$. Then, this SIP becomes
\begin{equation*}
\left\{
\begin{array}{cl}
	\min& (x_1+x_2-2)^2+(x_1-x_2)^2+30x^2_3 \\
	\st & -x_1w_1-x_2w_2+1\geq 0\quad\forall w \in W=\{w_1^2 + w_2^2 =1,\, w_1\geq 0\}, \\
	& x \in X=[100,100]^3\cap\{x_1 - x_2 - 2x_3 \ge 0, \,\,
            (x_1 - x_2 - 2x_3)^2 = (x_1 - x_2 )^2\} .
\end{array}
\right.
\end{equation*}
The minimizer for this equivalent GSIP is
$\tilde{x}\approx (0.7071, 0.7071, 0.0000)$,
so the minimizer for Example~\ref{Teo 1} is
$x^*\approx (0.7071, 0.7071)$.
\end{eg}

\begin{eg}\label{betro_sum}
The SIP from \cite{betro2004accelerated}:
\begin{equation*}
\left\{
\begin{array}{cl}
	\min & \sum\limits_{i=1}^{8} \frac{x_i}{i}\\
	\st & g(x,u) = \sum\limits_{i=1}^{8} u^{i-1}x_{i}-\frac{1}{2-u}\geq 0 \quad
	\forall u \in U=[0,1], \\
	&x \in X=[-1,1]^8 .
\end{array}
\right.
\end{equation*}
We replace $g(x,u)\geq 0$ by its equivalent form $(2-u)g(x,u)\geq 0$.
\end{eg}

\begin{eg}\label{Flouda 1} 	
The SIP for the Chebyshev approximation problem \cite{floudas2008adaptive,wang2015feasible}:
\begin{equation*}
\left\{
\begin{array}{cl}
	\min& x_{4}\\
	\st & \max \left|\sin (\pi u)-\left(x_{3} u^{2}+x_{2} u+x_{1}\right)\right|
	\leq x_4\quad\forall u\in U,\\
	&x\in X=\{x\in\mathbb{R}^4:
	-1\leq x_1\leq1, 3\leq x_2\leq 5, -5\leq x_3\leq -3, -1\leq x_4
	\leq 3\},\\
	& U=[0,1].
\end{array}
\right.
\end{equation*}
The first constraint is equivalent to
\[
x_4\pm\left(\sin (\pi u)-x_{3} u^{2}-x_{2} u-x_{1}\right) \geq 0,
\]
which gives a SIP. The function $\sin(\pi u)$
is approximated by the Taylor expansion of degree $11$.
If other degree Taylor expansions for $\sin(\pi u)$ are used,
we get the minimum values
$0.7147$, $0.4059$, $0.1039$, $0.0296$, $0.0292$,
when the expansion degrees are
$1$, $3$, $5$, $7$, $9$ respectively.
\end{eg}

\begin{eg}\label{Zakovic-minmax}
The SIP for the min-max problem \cite{zakovic2003semi}:
\begin{equation*}
\left\{
\begin{array}{cl}
	\min\limits_{x}\max\limits_{u}& 4(x_1-2)^2-2u^2_1+x^2_1u_1-u^2_2+2x^2_2u_2\\
	\st & x\in X=[-100,100]^2,\\
	&u\in U=[-5,5]^2.
\end{array}
\right.
\end{equation*}
Introduce a new variable $x_3$ such that the problem convert to an equivalent SIP problem:
\begin{equation*}
\left\{
\begin{array}{cl}
	\min& x_3\\
	\st &4(x_1-2)^2-2u^2_1+x^2_1u_1-u^2_2+2x^2_2u_2\leq x_3 \\
	&x\in X=[-100,100]^2,\\
	&u\in U=[-5,5]^2.
\end{array}
\right.
\end{equation*}
\end{eg}

\begin{eg}\label{lemondis_paper}
The SIP from \cite{bhattacharjee2005global}
\begin{equation*}
	\left\{
	\begin{array}{cl}
		\min & \frac{x^2_1}{3}+\frac{x_1}{2}+x^2_2-x_2\\
		\st & g(x,u)=\sin(u)-x^2_1-2x_1x_2u^2\geq 0\quad\forall u\in U,\\
		& x \in X=[-1000,1000]^2\\
		& U=[0,2].
	\end{array}
	\right.
\end{equation*}
The function $\sin(u)$ is approximated by the Taylor expansion of degree $7$.
If other degree Taylor expansions for $\sin(u)$ are used,
we get the same minimum value $-0.2500$,
when the expansion degrees are
$1$, $3$, $5$ respectively.
\end{eg}

\section{Some GSIP examples from references}
\label{AppendixB}

The GSIPs in this appendix are from existing references.
Their computational results are shown in
Table~\ref{all result for appendix B}.

\begin{eg}\label{Aboussoror1}
The GSIP from \cite{Aboussoror1}:
\begin{equation*}
	\left\{
	\begin{array}{cl}
		\min & x^2+|x|+1\\
		\st  & g(x,u)=x^2+2x+u_1+u_2\geq 0\quad\forall u\in U(x),\\
		& x \in X=[-\frac{1}{2},\frac{1}{2}],\\
		& U(x) \, =  \, \big \{u\in \mathbb{R}^2:  x^2-1\leq u_1 \leq -\frac{3}{4},\,
		x-1\leq u_2\leq -\frac{1}{2} \big \}.
	\end{array}
	\right.
\end{equation*}
Let $x_1:=x$, then introduce a new variable $x_2$ such that $x_2=x^2_1,\,x_2\geq 0$.
Let $x\in [-\frac{1}{2},\,\frac{1}{2}]^2\cap \{x_2=x^2_1,\,x_2\geq 0\}$,
then the objective function becomes $x_1+x_2+1$, and the constraints
$x^2_1-1\leq -\frac{3}{4},\,x_1-1\leq -\frac{1}{2}$ gives $U(x)$ nonempty with box constraints.
The polynomial extension is available as for box constraints.
The minimizer for this equivalent problem with new variable is
$\tilde{x}\approx ( 0.5000,    0.5000)$,
so the minimizer for Example~\ref{Aboussoror1} is $x^*\approx 0.5000$.
\end{eg}

\begin{eg}\label{Alexander 1}
The GSIP from \cite{articleAngelos2015,ruckmann2001second}:
\begin{equation*}
\left\{
\begin{array}{cl}
	\min & -0.5x^4_1+2x_1x_2-2x^2_1\\
	\st  & g(x,u)=-x^2_1+x_1+x_2-u^2_1-u^2_2\geq 0\quad\forall u\in U(x),\\
	& x \in X=[0,1]^2,\\
	& U(x) =\{u\in \mathbb{R}^3:  0\leq u\leq e, x_1-u^Tu\geq 0\}.
\end{array}
\right.
\end{equation*}
Make substitution $y_1\coloneqq u^2_1, y_2\coloneqq u^2_2, y_3\coloneqq u^2_3$,
then this GSIP becomes
\begin{equation*}
\left\{
\begin{array}{cl}
	\min & -0.5x^4_1+2x_1x_2-2x^2_1\\
	\st  & -x^2_1+x_1+x_2-y_1-y_2\geq 0\quad\forall y\in Y(x),\\
	& x \in X=[0,1]^2,\\
	& Y(x) = \{ y\in \mathbb{R}^3:  0\leq y \leq 1, e^Ty  \leq  x_1\}.
\end{array}
\right.
\end{equation*}
For all $x\in X$, we have $0\leq x_1\leq 1$, so the constraining set is nonempty and
$Y(x)=\{0 \leq y,e^Ty \leq x_1\},$ which has simplex constraints. The polynomial extension
is available as for simplex constraints. This new problem gives
$y^*\approx (  0.5000 ,   0.5000 ,   0.0000)$,
so Example~\ref{Alexander 1} gives $u^*\approx ( 0.7071,    0.7071 ,   0.0000)$.
\end{eg}

\begin{eg}\label{M.Diehl 1}
The GSIP from \cite{diehl2013lifting}:
\begin{equation*}
	\left\{
	\begin{array}{cl}
		\min  &  (x_1+\sqrt{2}+1)^2+(x_2-1)^2\\
		\st & g(x,u)=u\geq 0\quad\forall u\in U(x),\\
		& x\in X=\mathbb{R}^2,\\
		& U(x)=\{u \in \re :  x_1 \leq u \leq  x_2\}.
	\end{array}
	\right.
\end{equation*}
Clearly, when $x_1 > x_2$, the set $U(x)$ is empty.
We solve this GSIP in two cases.

\medskip
\noindent
Case I: If $x\in X_1 \coloneqq X\cap \{x_1\leq x_2\}$,
then $U(x)=\{x_1\leq u\leq x_2\}$ is a box constraint.
The polynomial extension is available as in \reff{PE:box}.
Solving this case of the GSIP by Algorithm~\ref{algorithm},
we get $x^*\approx ( 0.0000,    0.9999),\,f^*\approx  5.8284.$

\medskip
\noindent
Case II: If $x\in X_2 \coloneqq X \cap \{x_1> x_2\}$, then $U(x) = \emptyset $
and this GSIP is equivalent to minimizing $f(x)$ over the set
$X_2$. The strict inequality $x_1 > x_2$ is treated as $x_1  \ge  x_2$.
Solving this case of the polynomial optimization,
we get the optimal $\tilde{x}\approx(4.4,\, 5.3 )\cdot 10^{-14}$, $f(\tilde{x})\approx 6.8284$.
The inequality $x_1 \ge  x_2$ is active at $x^*$. This point $x^*$ is feasible for the GSIP,
but it is not optimal when compared with Case I.
\end{eg}

\begin{eg}\label{jongen 1}
The GSIP  from \cite{jongen1998generalized}:
\begin{equation*}
	\left\{
	\begin{array}{cl}
		\min& x_2\\
		\st & g(x,u)=u^3-x_2\geq 0\quad\forall u\in U(x),\\
		&x\in X=\mathbb{R}^2,\\
		&U(x)=\{u\in\mathbb{R} :  u\leq 0, -x^2_1+2x_2-u^3\leq 0\}.
	\end{array}
	\right.
\end{equation*}
Introduce the new variable $x_3$ such that $x^3_3=-x^2_1+2x_2$,
and add this constraint to constraining set $X$. Compare $x_3$ with $0$.
When $x_3\leq 0$, then $U(x) = \{x_3 \leq u \leq 0\}$ is a box constraint.
When $x_3 > 0$, the set $U(x)$ is empty.
We solve the GSIP in two cases.

\medskip
\noindent
Case I: For $x\in X_1 \coloneqq  X \cap \{x_3\leq 0, \,x^3_3=-x^2_1+2x_2\}$,
$U(x)$ is a nonempty box constraint.
The polynomial extension is available as in \reff{PE:box}.
Use Algorithm~\ref{algorithm} to solve this case of the GSIP, we get the solution
$(0.0000, \, 0.0000, \, -0.0037),\,f^*\approx -1.0\cdot 10^{-7}$.
So the minimizer for this case is $x^*=(0.0000, 0.0000).$

\medskip
\noindent Case II:
For $x\in  X_2 \coloneqq  X \cap \{x_3> 0, \,x^3_3=-x^2_1+2x_2\}$, $U(x)=\emptyset $.
The GSIP is the equivalent to minimize $f(x)$ over $X_2$.
Solving this case polynomial optimization, we get the minimizer
$(0.0000,\,0.0000,\,0.0000)$.
By comparison with Case I, the minimizer for this GSIP is $(0.0000, 0.0000).$
\end{eg}

\begin{eg}\label{Alexander_4}
The GSIP from \cite{articleAngelos2015}:
\begin{equation*}
\left\{
\begin{array}{cl}
	\min&x^2_1+x^2_2+x^2_3 \\
	\st \,& g(x,u)=-x_1(u_1+u^2_2+1)-x_2(u_1u_2-u^2_2)-x_3(u_1u_2+u^2_2+u_2)-1\geq 0\\
	&\forall u\in U(x) \coloneqq \{u\in\mathbb{R}^2 :
	0\leq u\leq 1,  x^2_1\leq u^2_1\}, \\
	& x\in X=[-1,0]^3 .
\end{array}
\right.
\end{equation*}
Since $x\in X=[-1,0]^3$, then the set $U(x)$ is equivalent to the box type constraint
$\{-x_1\leq u_1\leq 1, 0\leq u_2 \leq 1\}$.
The polynomial extension is available as for the box case.
\end{eg}

\begin{eg}\label{Alexander_5}
The GSIP from \cite{articleAngelos2015}:
\begin{equation*}
\left\{
\begin{array}{cl}
	\min&x^2_1+x^2_2+x^2_3 \\
	\st & g(x,u)= 2\sin(4u)-x_1-x_2\exp(x_3u)-\exp(x_2u)\geq 0\quad\forall u\in U(x),\\
	&x\in X=[-5, 5]^3,\\
	&U(x)=\{u\in\mathbb{R}:  0\leq u\leq 1, x_2+1-2u\geq 0\}.
\end{array}
\right.
\end{equation*}
We approximate $\sin(4u)$ by its Taylor expansion of degree $5$
and approximate $\exp(x_i u)$ by the Taylor expansion of degree $4$.
The set $U(x) = \emptyset$ when $\frac{1}{2}(x_2+1)> 1$.
Comparing $\frac{1}{2}(x_2+1)$ with $0$ and $1$ respectively,
we solve this GSIP in three cases.

\smallskip
\noindent
Case I: If $x \in X_1 \coloneqq  X\cap\{-1\leq x_2\leq 1\}$, then
$U(x)=\{0\leq u \leq \frac{1}{2}(x_2+1)\}$ is a box constraint.
Solving this case of the GSIP with Algorithm~\ref{algorithm}, we get
\[
x^*\approx ( -0.5000,   -0.5000,    0.0000),\quad f^*\approx 0.5000.
\]

\smallskip
\noindent
Case II: If $x\in X_2 \coloneqq X\cap \{x_2>1\}$, then $U(x)=[0,\,1]$ is a constant set
and the GSIP reduces to an SIP. The constraint $x_2>1$ is treated as $x_2\geq 1$
in the polynomial optimization. Solving this case by Algorithm~\ref{algorithm},
we get $\tilde{x}\approx(-4.2248,1.0000,-0.8737),\,f(\tilde{x})\approx  0.8686$.
This solution is not optimal when compared with Case I.

\smallskip
\noindent
Case III: If $x\in  X_3 \coloneqq X \cap \{x_2< -1\}$, then $U(x)=	\emptyset $
and the GSIP is reduced to minimizing $f(x)$ over the set $X_3$.
The constraint $x_2 < -1$ is treated as $x_2 \leq -1$.
Solving this polynomial optimization gives
$\bar{x}\approx(  0.0000,   -1.0000,    0.0000)$, $f(\bar{x})\approx  1.0000$,
which is not optimal when compared with Case I.

\smallskip
\noindent
If other degree Taylor expansions for $\sin(4u)$ are used,
we get the same minimum value $0.5000$
when the expansion degrees are
$1$, $3$, $5$ respectively
(the expansion degree for $\exp(x_i u)$ is one lower than that for $\sin(4u)$).
\end{eg}

\begin{eg}\label{Alexander_6}
The GSIP from \cite{articleAngelos2015}:
\begin{equation*}
\left\{
\begin{array}{cl}
\min & x_1+x_2 \\
\st & g(x,u)=u\geq 0\quad\forall u\in U(x),\\
 & x\in X=[-1, 1]^2,\\
 & U(x)=\{u\in\mathbb{R}:  -1\leq u\leq 1, x_1\leq u, x_2\leq u\}.
\end{array}
\right.
\end{equation*}
\noindent Introduce the new variable $x_3 \coloneqq \max(-1, x_1, x_2)$, which is equivalent to
\[
(x_3+1)(x_3-x_1)(x_3-x_2) = 0, \quad x_3+1 \ge 0, \quad x_3-x_1 \ge 0, \quad x_3 - x_2 \ge 0.
\]
We add the above constraints to description of the set $X$,
then $U(x)$ is equivalently given by the box constraint
$x_3 \leq u  \leq 1$. For $x_1, x_2 \in [-1,1]^2$,
we have $x_3\in [-1,1]$ and $U(x)\not=\emptyset$.
Solving the polynomial optimization relaxation,
we get two minimizers:
\[
(-1.0000,0.0000,0.0000), \quad (0.0000,-1.0000, 0.0000)
\]
with the same minimal value $f^*= -1.0000$.
So the minimizers for this GSIP are
\[
(-1.0000,0.0000), \quad (0.0000, -1.0000).
\]
\end{eg}

\begin{eg}\label{Alexander_7}
Consider the following GSIP from \cite{lemonidis2008global,articleAngelos2015}
\begin{equation*}
\left\{
\begin{array}{cl}
	\min&\exp(x_1)+\exp(x_2)+\exp(x_3) \\
	\st & g(x,u)=x_1+x_2u+x_3u^2-\frac{1}{1+u^2}\geq 0\quad\forall u\in U(x),\\
	&x\in X=[-1, 1]^3,\\
	&U(x)=\{u\in\mathbb{R}: 0  \leq u \leq 1, \,  2(x_2+x_3)\leq u \}.
\end{array}
\right.
\end{equation*}
We replace $\exp(x_i)$ by its $4$th-degree Taylor expansion
and replace $g(x,u)\geq 0$ by its equivalent form $(1+u^2)g(x,u)\geq 0$.
Comparing $2(x_2+x_3)$ with $0$ and $1$,
we solve this GSIP in three cases:

\smallskip
\noindent
Case I: If $x\in X_1 \coloneqq X\cap \{0  \le  x_2+x_3 \le  \frac{1}{2}\}$,
then $U(x)=\{2(x_2+x_3)\leq u\leq 1\}$ is a box constraint.
Applying Algorithm~\ref{algorithm} for the set $X_1$,
we get the solution
\[
\tilde{x}\approx (0.0000,\,0.2500,\,0.2500) ,\quad f(\tilde{x}) \approx 3.5680 .
\]
It is feasible but not optimal when compared with Case III.

\medskip
\noindent
Case II: If $x\in X_2 \coloneqq X\cap \{x_2+x_3< 0\}$,
then $U(x)=[0,\,1]$ and the GSIP reduces to an SIP.
Applying Algorithm~\ref{algorithm} for the set $X_2$, we get the solution
\[
 \bar{x}\approx (1.0000,\, -0.0002,\, -0.4998),\quad f(\bar{x})\approx 4.3246 .
\]
It is feasible but not optimal when compared with Case III.

\smallskip
\noindent
Case III: If $x\in  X_3 \coloneqq X \cap \{x_2+x_3> \frac{1}{2}\}$,
then $U(x)=\emptyset$ and the GSIP is equivalent to
minimizer $f(x)$ over the set $X_3$.
The constraint $x_2+x_3 > \frac{1}{2}$ is treated as
$x_2+x_3\geq  \frac{1}{2}$.  Solving this case polynomial optimization,
we get $x^*\approx (  -1.0000,    0.2500,    0.2500)$
and $f^*\approx 2.9361$.  By comparison of objective values,
this point is the minimizer for the GSIP.

\noindent
If other degree Taylor expansions for $\exp$ are used,
we get the minimum values
$2.5000$, $3.0625$, $2.9010$, $2.9340$, $2.9347$, $2.9361$
when the expansion degrees are
$1$, $2$, $3$, $4$, $5$, $6$ respectively.
\end{eg}

\begin{eg}\label{Harwood 1}
The GSIP from \cite{harwood2016lower}:
\begin{equation*}
\left\{   \begin{array}{cl}
	\min\limits&-(x_3-x_1)(x_4-x_2) \\
	\st & g(x,u)=\cos(u_1)\sin(u_2)- \frac{u_1}{u_2^2+1} -
	\sum\limits_{i=1}^2 \alpha_{i}u_{i}(1-u_{i}) \geq 0 \quad \forall u\in U(x),\\
	&x\in X,
\end{array}
\right.
\end{equation*}
where $\alpha_1=1.841, \alpha_2=6.841$.
The constraining sets are:
\[
\begin{array}{ll}
	X=\{x\in\mathbb{R}^4 :  u+e\geq 0, e-u\geq 0,x_1\leq x_3, x_2\leq x_4\},\\
	U(x)=\{ u \in\mathbb{R}^2:  x_1\leq u_1\leq x_3, x_2\leq u_2\leq x_4\}.
\end{array}
\]
We use the formula
$\cos(u_1)\sin(u_2)=\frac{1}{2}[\sin(u_2-u_1)+\sin(u_1+u_2)]$
and approximate the $\sin$ function by its $5$th-degree Taylor expansion.
The inequality $g(x,u)\geq 0$ is replaced
by its equivalent form $(u^2_2+1)g(x,u)\geq 0$.
For all $x\in X$, we have $x_1\leq x_3, \, x_2\leq x_4$,
so the constraining set is nonempty and has box case constraint.
The polynomial extension is available
as for the box type constraint.
If other degree Taylor expansions for $\sin(u_2-u_1)$, $\sin(u_2+u_1)$ are used,
we get the minimum values
$-1.6880$, $-1.6985$, $-1.6980$, $-1.6980$
when the expansion degrees are
$1$, $3$, $5$, $7$ respectively.
\end{eg}

\begin{eg}\label{O.Stein 1}
The GSIP from \cite{stein2003solvingpaper}:
\begin{equation*}
\left\{\begin{array}{cl}
	\min&-\pi x_3x_4 \\
	\st &
	\begin{bmatrix}
		u_1+u^2_2\\
		3-u_1-4u_2\\
		u_2+1\\
	\end{bmatrix}
	\geq 0\quad \forall u\in U(x),\\
	&x\in X,
\end{array} \right.
\end{equation*}
where the constraining sets
\[
\begin{array}{rcl}
	X &=& \left\{	x\in\mathbb{R}^4 :
	\begin{array}{ll}
		100-x^Tx\geq 0, x_3\geq 0, x_4\geq 0
	\end{array}\right\},\\
	U(x) &=& \left\{u\in\mathbb{R}^2\,\middle\vert\,\begin{array}{ll}
		\frac{(u_1-x_1)^2}{x^2_3}+\frac{(u_2-x_2)^2}{x^2_4} \leq 1
\end{array}\right\}.
\end{array}
\]
The polynomial extension is available
as for the ellipsoidal type constraint.
\end{eg}

\section{A numerical certificate for Example~\ref{exmp:4.5} }
\label{AppendixC}

Note that $c_1(x) = x_2-x_1$, $c_2(x) = x_1$, $c_3(x) = 2-x_1-x_2$,
$h_1(x,u) = x_2-x_1u_1-x_2u_2$, $h_2(x,u) = u_1+x_1u_2-x_1x_2$,
$h_3(x,u) = x_1+3u_1-u_2$, $h_4(x,u) =u_1-x^2_1u_2$.
The polynomial extension $q = (\phi_1, \phi_2)$
given in Example~\ref{exmp:4.5}
satisfies \reff{polyext:q:num},
which is shown by the following representations:
\begin{align*}
& h_1(x,q(x))  \approx  \\
& \begin{bmatrix}
	1 \\ x_1 \\ x_2
\end{bmatrix}^T \left(\begin{bmatrix}
	-0.0003 & -0.0002 & -0.0001 \\
	-0.0002 & 0.3222 & -0.3224 \\
	-0.0001 & -0.3224 & 0.3220
\end{bmatrix} \cdot 10^{-6} + \begin{bmatrix}
	0.2898 & 0.1010 & -0.1109 \\
	0.1010 & 0.0352 & -0.0386 \\
	-0.1109 & -0.0386 & 0.0424
\end{bmatrix}\cdot c_1(x) \right. \\
& \left. + \begin{bmatrix}
	0.0193 & 0.0340 & -0.0533 \\
	0.0340 & 0.0601 & -0.0941 \\
	-0.0533 & -0.0941 & 0.1474
\end{bmatrix}\cdot c_2(x) + \begin{bmatrix}
	0.0000 & 0.0000 & 0.0000 \\
	0.0000 & 0.2939 & 0.0208 \\
	0.0000 & 0.0208 & 0.0424
\end{bmatrix}\cdot c_3(x) \right)\begin{bmatrix}
	1 \\ x_1 \\ x_2
\end{bmatrix},
\end{align*}

\begin{align*}
& h_2(x,q(x))  \approx  \\
& \begin{bmatrix}
	1 \\ x_1 \\ x_2
\end{bmatrix}^T \left(\begin{bmatrix}
	0.0001 & 0.0028 & -0.0029 \\
	0.0028 & 0.0938 & -0.0965 \\
	-0.0029 & -0.0965 & 0.0994
\end{bmatrix}  + \begin{bmatrix}
	0.0002 & 0.0005 & -0.0009 \\
	0.0005 & 0.0012 & -0.0024 \\
	-0.0009 & -0.0024 & 0.0047
\end{bmatrix}\cdot c_1(x)
\right. \\
& \left. + \begin{bmatrix}
	0.1229 & -0.0586 & -0.0643 \\
	-0.0586 & 0.0332 & 0.0253 \\
	-0.0643 & 0.0253 & 0.0389
\end{bmatrix}\cdot c_2(x)
+\begin{bmatrix}
	0.1352 & 0.0658 & 0.0253 \\
	0.0658 & 0.0320 & 0.0123 \\
	0.0253 & 0.0123 & 0.0047
\end{bmatrix}\cdot c_3(x)
\right)\begin{bmatrix}
	1 \\ x_1 \\ x_2
\end{bmatrix},
\end{align*}

\begin{align*}
&  h_3(x,q(x))   \approx  \\
& \begin{bmatrix}
	1 \\ x_1 \\ x_2
\end{bmatrix}^T \left(\begin{bmatrix}
	0.3031 & -0.0744 & -0.2719 \\
	-0.0744 & 0.9268 & 0.3227 \\
	-0.2719 & 0.3227 & 0.5376
\end{bmatrix}\cdot 10^{-5}  + \begin{bmatrix}
	0.0151 & -0.0098 & -0.0300 \\
	-0.0098 & 0.0063 & 0.0194 \\
	-0.0300 & 0.0194 & 0.0596
\end{bmatrix}\cdot c_1(x)
\right. \\
& \left. + \begin{bmatrix}
	0.3705 & 0.2022 & 0.2995 \\
	0.2022 & 0.1215 & 0.1567 \\
	0.2995 & 0.1567 & 0.2463
\end{bmatrix}\cdot c_2(x)
+\begin{bmatrix}
	0.0506 & -0.0764 & -0.0550 \\
	-0.0764 & 0.1152 & 0.0829 \\
	-0.0550 & 0.0829 & 0.0596
\end{bmatrix}\cdot c_3(x)
\right)\begin{bmatrix}
	1 \\ x_1 \\ x_2
\end{bmatrix},
\end{align*}

\begin{align*}
& h_4(x,q(x))   \approx  \\
& \begin{bmatrix}
	1 \\ x_1 \\ x_2
\end{bmatrix}^T \left(\begin{bmatrix}
	0.0517 & -0.0407 & -0.0222 \\
	-0.0407 & 0.0320 & 0.0175 \\
	-0.0222 & 0.0175 & 0.0095
\end{bmatrix}
+ \begin{bmatrix}
	0.1551 & 0.0706 & 0.0123 \\
	0.0706 & 0.0322 & 0.0056 \\
	0.0123 & 0.0056 & 0.0010
\end{bmatrix}
\cdot c_1(x)
\right. \\
& \left. + \begin{bmatrix}
	0.1930 & -0.2930 & 0.0303 \\
	-0.2930 & 0.4448 & -0.0461 \\
	0.0303 & -0.0461 & 0.0048
\end{bmatrix}
\cdot c_2(x)
+\begin{bmatrix}
	0.1094 & -0.0753 & -0.0102 \\
	-0.0753 & 0.0518 & 0.0070 \\
	-0.0102 & 0.0070 & 0.0010
\end{bmatrix}
\cdot c_3(x)
\right)\begin{bmatrix}
	1 \\ x_1 \\ x_2
\end{bmatrix}.\\
\end{align*}
The above matrices all positive semidefinite (up to numerical errors).
For neatness, only four decimal digits are shown.

\end{document}